\documentclass[10pt]{article}
\usepackage{mystyle}
\usepackage{tikz}
\usepackage{tikz-network}
\usetikzlibrary{matrix,arrows,positioning,calc,shapes}

\DeclareMathOperator{\cl}{cl}
\newcommand{\sPG}{\mathsf{PG}}
\newcommand{\curve}{\mathbf{r}}  % set some notation for a parameterized curve which we can agree on and change later
\newcommand{\library}{{NDMA}} 

\newcommand{\concat}{\textbf{V}}

\providecommand{\keywords}[1]
{
  \small    
  \textit{\textit{Keywords:}} Network dynamical systems, Bifurcation analysis, High dimensional ODE, Numerical analysis of dynamical systems
}
\providecommand{\subjclass}[1]
{
  \small    
  \textit{\textit{2020 MSC:}} 37M20,  37G35, 65L99
}

\title{Global analysis of regulatory network dynamics: equilibria and saddle-node bifurcations}
\author{Shane Kepley\thanks{Department of Mathematics, VU Amsterdam, 1081 HV Amsterdam, The Netherlands (s.kepley@vu.nl).},  Konstantin Mischaikow\thanks{Department of Mathematics, Rutgers University, Piscataway, NJ 08854 USA (mischaik@math.rutgers.edu)}, and Elena Queirolo\thanks{IRMAR, Université de Rennes, 35000 Rennes, France  (elena.queirolo@univ-rennes.fr)}}
\date{\today}

\begin{document}

\maketitle

\begin{abstract}
In this paper we describe a combined combinatorial/numerical approach to studying equilibria and bifurcations in network models arising in Systems Biology. ODE models of the dynamics suffer from high dimensional parameters which presents a significant obstruction to studying the global dynamics via numerical methods. The main point of this paper is to demonstrate that adapting and combining classical techniques with recently developed combinatorial methods  provides a richer picture of the global dynamics despite the high parameter dimension. 

Given a network topology describing state variables which regulate one another via monotone and bounded functions, we first use the {\em Dynamic Signatures Generated by Regulatory Networks} (DSGRN) software to obtain a combinatorial summary of the dynamics. This summary is coarse but global and we use this information as a first pass to identify ``interesting'' subsets of parameters in which to focus. We construct an associated ODE model with high parameter dimension using our {\em Network Dynamics Modeling and Analysis} (NDMA) Python library. We introduce algorithms for efficiently investigating the dynamics in these ODE models restricted to these parameter subsets. Finally, we perform a statistical validation of the method and several interesting dynamical applications including finding saddle-node bifurcations in a $54$ parameter model. 

\end{abstract}

\section{Introduction}
\label{sec:introduction}

Consider an explicit differential equation
\begin{equation}
    \label{eq:abstract}
    \dot{x} = f(x,\lambda), \quad x\in \rr^N,\ \lambda \in \rr^M,
\end{equation}
where $x$ is a vector of state variables and $\lambda$ is a vector of parameters. The most ambitious goal is to understand the {\em global dynamics} of this system, by which we mean a complete understanding of the possible dynamics over all parameter values. In all but the simplest cases this goal is infeasible. A more modest goal is to choose a relevant invariant set, e.g., an equilibrium or periodic orbit, and identify the set of parameters at which this exists. This leads to the closely related problem of identifying bifurcations, i.e., sets of parameters at which the topology or dynamics of this invariant set changes. Bifurcation analysis is an active area of research in dynamical systems theory and we do not intend to present an in depth survey of the enormous amount of literature on the subject. An introduction to the theory can be found in \cite{chow:hale, wig90, Kuznetsov1998, palis:takens} and a thorough treatment of modern computational methods for many bifurcations in \cite{Dankowicz}. 
Compared to these works the focus of this paper is much narrower; we address only the question of existence and stability of equilibrium solutions and their saddle-node bifurcations in a class of ordinary differential equations (ODEs) modeling gene regulatory network dynamics.  
However, our approach is fundamentally different. 

Our motivation for developing an alternative approach comes from the need to analyze systems for which $M$, the dimension of parameter space, is large, and the relevant parameter values are unknown.
To put this into context we note that if $M=1$, then the existence of equilibria is often approached by identifying an equilibrium at a given parameter value and then using continuation techniques to characterize the equilibria as a function of $\lambda$, i.e., a parameterized curve in $\rr^N$. 
Local bifurcations can be identified as isolated points in parameter space where the equilibria disappear or change stability. 
Despite its simplicity, this approach is extremely powerful and has been widely applied to studying codimension-one bifurcations. 

This same strategy can, at least in principle, be carried out when $M \geq 2$. 
However, this presents practical and computational challenges \cite{gameiro:lessard:pugliese, henderson}. 
For instance, if $M=2$ then characterizing the equilibria is already more complicated as one needs to identify a two-dimensional surface of points instead of a curve. 
Similarly, identifying bifurcations of these equilibria also becomes increasingly difficult. 
In this context saddle-node bifurcations (associated with the appearance and disappearance of equilibria) generically take the form of $1$-dimensional curves embedded in this surface that themselves can be identified via continuation. 
We note that already this is much harder than for $M=1$.

In contrast, consider Figure~\ref{fig:EMT_model} that represents a simplified version of a regulatory network for Epithelial-Mesenchymal Transition (EMT) that is of oncological interest \cite{Xin2020MultistabilityIT}.
EMT acts as a switch, the outcome of which determines cellular phenotypes.
The simplest classic robust switch involves hysteresis, and thus, is associated with changes in the number of stable equilibria and the existence of saddle-node bifurcations.
As is explained in Section~\ref{sec:hill models} it is reasonable to assume that an ODE representing this EMT model is a six-dimensional system with a  $54$-dimensional parameter space.
Furthermore, at best the relevant parameter values are known imprecisely.
Thus, there is no natural starting point for finding equilibria, nor are there natural directions in which to perform continuation.

\begin{figure}
	\centering
\resizebox{0.5\textwidth}{!}{ % resize here as needed
\begin{tikzpicture} 
% fake vertices to attach edges should be slightly larger
\Vertex[x = -1, y = 1, size = 0.9, opacity = 0, Pseudo = True]{x1} 
\Vertex[x = 0, y = 0, size = 0.9, opacity = 0, Pseudo = True]{x2} 
\Vertex[x = 1, y = 1, size = 0.9, opacity = 0, Pseudo = True]{x3} 
\Vertex[x = -1.5, y = -0.5, size = 0.9, opacity = 0, Pseudo = True]{x4} 
\Vertex[x = 0, y = -1.5, size = 0.9, opacity = 0, Pseudo = True]{x5} 
\Vertex[x = 2, y = -0.5, size = 0.9, opacity = 0, Pseudo = True]{x6} 
% real vertices with labels 
\Vertex[x = -1, y = 1, size = 0.8, opacity = 0, label = $1$]{X1} 
\Vertex[x = 0, y = 0, size = 0.8, opacity = 0, label = $2$]{X2} 
\Vertex[x = 1, y = 1, size = 0.8, opacity = 0, label = $3$]{X3} 
\Vertex[x = -1.5, y = -0.5, size = 0.8, opacity = 0, label = $4$]{X4} 
\Vertex[x = 0, y = -1.5, size = 0.8, opacity = 0, label = $5$]{X5} 
\Vertex[x = 2, y = -0.5, size = 0.8, opacity = 0, label = $6$]{X6} 
% edges
\Edge[lw = 1, color = black, Direct, bend = 15](x1)(x3)
\Edge[lw = 1, color = black, Direct, bend = 15, style = {-{|[scale = 0.5]}}](x2)(x1)
\Edge[lw = 1, color = black, Direct, bend = 15, style = {-{|[scale = 0.5]}}](x2)(x5)
\Edge[lw = 1, color = black, Direct, bend = 15, style = {-{|[scale = 0.5]}}](x3)(x2)
\Edge[lw = 1, color = black, Direct, bend = 15](x3)(x5)
\Edge[lw = 1, color = black, Direct, bend = 15, style = {-{|[scale = 0.5]}}](x3)(x6)
\Edge[lw = 1, color = black, Direct, bend = 15, style = {-{|[scale = 0.5]}}](x4)(x1)
\Edge[lw = 1, color = black, Direct, bend = 15, style = {-{|[scale = 0.5]}}](x4)(x5)
\Edge[lw = 1, color = black, Direct, bend = 15, style = {-{|[scale = 0.5]}}](x5)(x2)
\Edge[lw = 1, color = black, Direct, bend = 15, style = {-{|[scale = 0.5]}}](x5)(x4)
\Edge[lw = 1, color = black, Direct, bend = -15, style = {-{|[scale = 0.5]}}](x5)(x6)
\Edge[lw = 1, color = black, Direct, bend = 15, style = {-{|[scale = 0.5]}}](x6)(x3)
\end{tikzpicture} 
}
	\caption{The network structure for the EMT model.}
\label{fig:EMT_model}
\end{figure}

This motivates us to formulate the  {\em equilibrium} problem  
\begin{description}
    \item[EQ:] Reliably count and approximate all equilibrium solutions which occur for a given parameter value.
\end{description}
and the {\em saddle-node bifurcation} problem
\begin{description}
    \item[SNB:] Reliably predict whether a given parameter is ``near'' an $M-1$ dimensional saddle-node bifurcation surface and if so, obtaining a nearby parameter which lies on this surface.
\end{description}
Observe that in the context of EMT this information tells us parameter values at which single or multiple phenotypes may occur and parameter values at which dramatic changes in phenotype may occur.

The main idea of our approach is to consider Equation \eqref{eq:abstract} at two ``resolutions.'' 
\begin{description}
    \item[Step 1] We apply combinatorial methods that have been implemented in the Dynamic Signatures Generated by Regulatory Networks (DSGRN) library \cite{DSGRN}. 
For the purposes of this paper DSGRN takes as input a regulatory network, e.g., Figure~\ref{fig:EMT_model}, and as output provides a decomposition of parameter space into explicit semi-algbraic sets and for each such semi-algbraic set a count indicating the expected number of stable equilibria.  
The hypersurfaces separating semi-algebraic sets with different numbers of expected stable fixed points are used to provide educated guesses of parameter values where saddle-node bifurcations occur.
\item[Step 2] We choose a class of ODEs and use the information gained from this combinatorial analysis to guide a refined numerical investigation of the dynamics using our Network Dynamics Modeling and Analysis (NDMA) software \cite{NDMA}.
NDMA is an open source Python library that has been developed to efficiently perform numerical analyses for the class of ODEs used in this work, but is applicable to a much broader family of equations.
\end{description}
 
The thesis of this paper is that by combining combinatorial and numerical techniques we obtain insight which neither method achieves alone. In particular, we demonstrate the ability to efficiently address both the equilibrium and saddle-node bifurcation problems described above for a remarkably large class of network models. 

The reader may be surprised by the fact that the ODE of interest can be chosen in {\bf Step 2}.
However, DSGRN has been designed  with the goal of characterizing dynamics for a broad class of dynamical systems as opposed to finding specific solutions to any particular system.
A systematic description of DSGRN goes beyond the scope of this paper (see \cite{cummins:gedeon:harker:mischaikow:mok, DSGRN, gedeon:cummins:harker:mischaikow, kepley:mischaikow:zhang, inpreparation:riveira:mishaikow}), and thus, in Section~\ref{sec:DSGRN} we limit ourselves to a brief description using a particularly simple network commonly referred to as the \emph{Toggle Switch}.

To quantify the usefulness of the DSGRN output requires that we consider a particular family of ODEs.
In \cite{inpreparation:riveira:mishaikow} it is shown that for a broad but carefully chosen family of ODEs called ramp systems, defined in terms of nonlinearities build from piecewise linear functions of a single variable, the DSGRN results can be identified with the dynamics of the ODE. 
For this paper we use Hill models (see Section~\ref{sec:hill models} for an explicit formulation).
There are three important reasons that we choose to work with these models.
\begin{enumerate}
    \item They are expressed in terms of rational functions (a more common modeling choice than ramp systems).  
    As a consequence, a numerical library can be built that computes derivatives of the vector field up to any given order, with respect to both the phase space variables and any of the parameters.
    Therefore, this library provides a basis for rigorous validation of the numerically detected dynamics.
    \item The DSGRN parameter space embeds in a very natural way into the Hill model parameter space.   
    \item Hill functions are commonly used in the context of systems biology.
\end{enumerate}

In Section~\ref{sec:equilibria} we turn to the NDMA computations and begin by providing a four step process for identifying equilibria and saddle-node bifurcations.
This includes defining the paths in parameter space along which we search for the bifurcations.
It is worth emphasizing that for reasons of computational efficacy we do not use continuation, but rather introduce an algorithm for obtaining \emph{a priori} bounds  on equilibria.
Both equilibria and saddle-node bifurcations are characterized as zero finding problems.
In the case of equilibria this allows us to use the so called ``radii polynomials'' to equate or distinguish numerical zeros.

In Section~\ref{sec:toggle_switch} we apply NDMA to the toggle switch where the associated ODEs are based on Hill function nonlinearities.
This is a two dimensional system of ODEs with a ten-dimensional parameter space.
We take advantage of the relative low dimensionality of this system  to try to explain in detail how DSGRN and NDMA can be used to extract important information about the dynamics. 
In fact, in hopes of clarity  we fix several parameters as this allows us to make the results more readable and easier to plot.
Since even for this reduced system it is impossible to fully explore the dynamics numerically we present some statistical information that is meant to justify the predictive ability of DSGRN while highlighting the subtlety of bifurcations even in this simple system.

In Section~\ref{sec:EMT_dyn} we turn our attention to the EMT network of Figure~\ref{fig:EMT_model} where the associated ODEs are again based on Hill function nonlinearities with the goal of demonstrating the capabilities of DSGRN and NDMA.
This is a six-dimensional ODE with a 54-dimensional parameter space.
We do not consider any form of dimension reduction as this is antithetical to the thesis of this work; our goal is to be able to systematically and efficiently identify bifurcations in high dimensional parameters spaces.
Again, we use statistical results to demonstrate that the DSGRN predictions can provide meaningful intuition with extremely low computational cost.

We conclude in Section~\ref{sec:conclusion} with a brief review of the main contributions of this paper.

\section{DSGRN and the Toggle Switch}
\label{sec:DSGRN}

As indicated in the introduction a systematic description of DSGRN exceeds the scope of this paper.
Thus, to provide intuition we focus on a particularly simple regulatory network commonly referred to as the Toggle Switch, but remark  that within the synthetic biology community the Toggle Switch represents a first successful design of a synthetic genetic switch \cite{collins} and therefore has significant biological relevance.

The input to DSGRN is a regulatory network, more precisely an annotated directed graph $(\cG,p)$ (directed edges take the form $\rightarrow$ indicating activation or $\dashv$ indicating repression) and a list of interaction functions (see Section~\ref{sec:hill models} for formal definitions).
For the Toggle Switch the regulatory network is shown in Figure~\ref{fig:TS}(a).

\begin{figure}[tbh]
\begin{center}
\fbox{\parbox{.95\textwidth}{
\begin{picture}(400,150)
\put(5,5){(a)}
\put(15,15){\begin{tikzpicture}
%			[main node/.style={circle,fill=white!20,draw,font=\sffamily\normalsize\bfseries},scale=2.5]
			[main node/.style={circle,fill=white!20,draw,font=\sffamily\tiny\bfseries},scale=2.5]
				\node[main node] (1) at (0,1) {1};
				\node[main node] (2) at (0,0) {2};

				\path[->,>=angle 90,thick]
				(1) edge[-|,shorten <= 2pt, shorten >= 2pt, bend left] node[] {} (2)
				(2) edge[-|,shorten <= 2pt, shorten >= 2pt, bend left] node[] {} (1)
				;
			\end{tikzpicture}
		}

\put(80,5){(b)}	
\put(90,15){\begin{tikzpicture}
		[main node/.style={circle,fill=white!20,draw,font=\sffamily\tiny\bfseries},scale=1.8]
				
		\node[main node] (1) at (0,2) {1};
		\node[main node] (2) at (1,2) {2};
		\node[main node] (3) at (2,2) {3};
		\node[main node] (4) at (0,1) {4};
		\node[main node] (5) at (1,1) {5};
		\node[main node] (6) at (2,1) {6};
		\node[main node] (7) at (0,0) {7};
		\node[main node] (8) at (1,0) {8};
		\node[main node] (9) at (2,0) {9};

		\path[->,>=angle 90,thick]
		(1) edge[-,shorten <= 2pt, shorten >= 2pt] node[] {} (2)
		(1) edge[-,shorten <= 2pt, shorten >= 2pt] node[] {} (4)
		(2) edge[-,shorten <= 2pt, shorten >= 2pt] node[] {} (3)
		(2) edge[-,shorten <= 2pt, shorten >= 2pt] node[] {} (5)
		(3) edge[-,shorten <= 2pt, shorten >= 2pt] node[] {} (6)
		(4) edge[-,shorten <= 2pt, shorten >= 2pt] node[] {} (5)
		(4) edge[-,shorten <= 2pt, shorten >= 2pt] node[] {} (7)
		(5) edge[-,shorten <= 2pt, shorten >= 2pt] node[] {} (6)
		(5) edge[-,shorten <= 2pt, shorten >= 2pt] node[] {} (8)
		(6) edge[-,shorten <= 2pt, shorten >= 2pt] node[] {} (9)
		(7) edge[-,shorten <= 2pt, shorten >= 2pt] node[] {} (8)
		(8) edge[-,shorten <= 2pt, shorten >= 2pt] node[] {} (9)
				;
			\end{tikzpicture}		
			}	

\put(267,65){(c)}
\put(277,75){
	\begin{tikzpicture}[scale=1]
	    \draw [->] (1,1) -- (0.1,0.1);
		\draw [->] (1,1) -- (1.9,0.1);
		
	    \draw [fill] (1,1) circle [radius=0.1];
		\node at  (1,1.3)  {$c$};
		\draw [fill] (0,0) circle [radius=0.1];
		\node at  (0.3,0)  {$a$};
		\draw [fill] (2,0) circle [radius=0.1];
		\node at  (1.7,0)  {$b$};
		;
	\end{tikzpicture}
}	

\put(240,5){(d)}
\put(250,15){
	\begin{tikzpicture}[scale=1]
		\draw [fill] (0,0) circle [radius=0.1];
		\node at  (0.3,0)  {$a$};
		
		;
	\end{tikzpicture}
}	

\put(295,5){(e)}
\put(305,15){
	\begin{tikzpicture}[scale=1]
		\draw [fill] (0,0) circle [radius=0.1];
		\node at  (0.3,0)  {$b$};
		
		;
	\end{tikzpicture}
	}
\put(350,5){(f)}
\put(360,18){
	\begin{tikzpicture}[scale=1]
		\draw [fill] (0,0) circle [radius=0.1];
%		\node at  (0.3,0)  {$b$};
		
		;
	\end{tikzpicture}	
}
\end{picture}
}}
\caption{(a) Toggle Switch. (b) Parameter graph for Toggle Switch. (c) Morse graph at node 5 of parameter graph. (d) Morse graph at nodes 1, 2, and 4 of parameter graph. Note that the element $a$ of these Morse graph can be directly identified with the element $a$ of the Morse graph at node 5. 
(e) Morse graph at nodes 6, 8, and 9 of parameter graph.
Note that the element $b$ of these Morse graph can be direcctly identified with the element $b$ of the Morse graph at node 5.
(f) Morse graph at nodes 3 and 7.
}
\end{center}
\label{fig:TS}
\end{figure}

DSGRN constructs a parameter space $\Xi =(0,\infty)^{N+3E}$, where $N$ and $E$ are the number of nodes and edges, respectively, in the directed graph.
Thus, for the Toggle Switch, the parameter space is $(0,\infty)^8$.

The output of DSGRN is the DSGRN database. 
This is best understood conceptually as an annotated graph that we refer to as a \emph{parameter graph} and denote by $\sPG(\cG,p)$.
Each node $v$ of $\sPG(\cG,p)$ corresponds to a semi-algebraic set $R(v)\subset\Xi =  (0,\infty)^{N+3E}$ given in terms of explicit inequalities, and the union of closures of all these semi-algebraic sets is $\overline{\Xi} = [0,\infty)^{N+3E}$.
If two nodes $v$ and $v'$ share an edge in the parameter graph, then their expressions as semi-algebraic sets differs by exactly one inequality, and in particular $\cl(R(v))\cap \cl(R(v'))$ is a codimension 1 hypersurface in $ [0,\infty)^{N+3E}$.
DSGRN decomposes the parameter space of the Toggle Switch into nine regions and the parameter graph for the Toggle Switch is shown in Figure~\ref{fig:TS}(b).

In addition, the DSGRN database assigns to each node $v$ a partially ordered set $(\sM(v),\leq_v)$, that is visually expressed via the Hasse diagram and called the \emph{Morse graph} at parameter node $v$.
The interpretation of the Morse graph in the context of continuous dynamics of an ODE (we emphasize that no ODE has been introduced yet) is that
the elements of the Morse graph represent potential recurrent dynamics and the partial order indicates the direction of the dynamics, i.e., if $p,q\in\sM$ and $p<q$, then no orbit can limit in backward time to the recurrent dynamics associated to $p$ and in forward time to the recurrent dynamics associated to $q$. 

Returning to the Toggle Switch, Figure~\ref{fig:TS}(c)-(f) indicate the Morse graphs at different nodes of the parameter graph.
In the context of an ODE an appropriate interpretation of this information is that at the parameter values associated with Node 5 (we denote this as $R(5)$), the system exhibits bistability where the elements $a$ and $b$ indicate the existence of distinct attractors and the element $c$ indicates the existence of an invariant set that acts as a separatrix between the two attractors.
At the remaining parameter nodes, the Morse graph contains a single element suggesting that at the associated parameter values the dynamics is monostable.

Ignoring for the moment the meaning of the parameters, the third column of Table~\ref{tab:parameter_regions} provides the explicit decomposition of parameter space.
There is an edge in the parameter graph between Node 5 and Node 2 and from Table~\ref{tab:parameter_regions} the reader can check that the shared boundary between the associated regions $R(5)$ and $R(2)$ is the $7$ dimensional surface defined by
\begin{equation}
    \label{eq:52face}
\ell_{1,2} = \gamma_1\theta_{2,1}, \quad  \ell_{2,1}  <\gamma_2\theta_{1,2}   < \ell_{2,1} + \delta_{2,1}.
\end{equation}
If the recurrent dynamics captured by the elements $a$, $b$, and $c$ of the Morse graph at Node 5 are equilibria, then the change in dynamics between parameters in Node 5 and Node 2 is most easily associated with the occurrence of a saddle-node bifurcation, i.e., for each parameter value that lies on the hypersurface defined by  \eqref{eq:52face} there exists a saddle-node bifurcation. 

Conceptually the DSGRN analysis of the EMT network of Figure~\ref{fig:EMT_model} is the same as that of the Toggle Switch.
However, the DSGRN parameter space for EMT is $(0,\infty)^{42}$ and the associated parameter graph has exactly 10,368,000,000 nodes each of which is explicitly known.
We remark that the DSGRN computations are sufficiently efficient that it is feasible to compute the Morse graph for each parameter value.
This efficacy is what allows us to carry out the experiments described in Section~\ref{sec:EMT_dyn}.

\begin{table*}\centering
	\begin{tabular}{@{}clcll@{}}
		\toprule
		Node & Phenotype & Attractors & Inequalities & Reduced Inequalities \\
		\midrule 
		% Region 1
		\multirow{2}{*}{$1$} &\multirow{2}{*}{Monostable } & \multirow{2}{*}{$\setof*{(1, 0)}$} & 
		$\gamma_1 \theta_{2,1} < \ell_{1,2} < \ell_{1,2} + \delta_{1,2}$ & $1 < \ell_{1,2}$ \\
		{} & {} & {} & 
		$\ell_{2,1}  < \ell_{2,1} + \delta_{2,1} < \gamma_2 \theta_{1,2}$ &
		$\ell_{2,1} + \delta_{2,1} < \gamma_2$ \\
		\midrule
		% Region 2
		\multirow{2}{*}{$2$} &\multirow{2}{*}{Monostable} & \multirow{2}{*}{$\setof*{(1, 0)}$} & 
		$\gamma_1 \theta_{2,1} < \ell_{1,2} < \ell_{1,2} + \delta_{1,2}$ & $1 < \ell_{1,2}$ \\
		{} & {} & {} & 
		$\ell_{2,1} < \gamma_2 \theta_{1,2} < \ell_{2,1} + \delta_{2,1}$ &
		$\ell_{2,1} < \gamma_2 < \ell_{2,1} + \delta_{2,1}$ \\
		\midrule
		% Region 3
		\multirow{2}{*}{$3$} &\multirow{2}{*}{Monostable} & \multirow{2}{*}{$\setof*{(1, 1)}$} & 
		$\gamma_1 \theta_{2,1} < \ell_{1,2} < \ell_{1,2} + \delta_{1,2}$ & $1 < \ell_{1,2}$ \\
		{} & {} & {} & 
		$\gamma_2 \theta_{1,2} < \ell_{2,1} < \ell_{2,1} + \delta_{2,1}$ & $\gamma_2 < \ell_{2,1}$ \\
		\midrule
		% Region 4
		\multirow{2}{*}{$4$} &\multirow{2}{*}{Monostable} & \multirow{2}{*}{$\setof*{(1, 0)}$} & 
		$\ell_{1,2} < \gamma_1 \theta_{2,1} < \ell_{1,2} + \delta_{1,2}$ & 
		$\ell_{1,2} < 1 < \ell_{1,2} + \delta_{1,2}$ \\
		{} & {} & {} & 
		$\ell_{2,1} < \ell_{2,1} + \delta_{2,1} < \gamma_2 \theta_{1,2}$ &
		$\ell_{2,1} + \delta_{2,1} < \gamma_2$ \\
		\midrule
		% Region 5
		\multirow{2}{*}{$5$} &\multirow{2}{*}{Bistable} & \multirow{2}{*}{$\setof*{(0,1), (1, 0)}$} & 
		$\ell_{1,2} < \gamma_1 \theta_{2,1} < \ell_{1,2} + \delta_{1,2}$ & 
		$\ell_{1,2} < 1 < \ell_{1,2} + \delta_{1,2}$ \\
		{} & {} & {} & 
		$\ell_{2,1} < \gamma_2 \theta_{1,2} < \ell_{2,1} + \delta_{2,1}$ &
		$\ell_{2,1} < \gamma_2 < \ell_{2,1} + \delta_{2,1}$ \\\midrule
		% Region 6
		\multirow{2}{*}{$6$} &\multirow{2}{*}{Monostable} & \multirow{2}{*}{$\setof*{(0, 1)}$} & 
		$\ell_{1,2} < \gamma_1 \theta_{2,1} < \ell_{1,2} + \delta_{1,2}$ & 
		$\ell_{1,2} < 1 < \ell_{1,2} + \delta_{1,2}$ \\
		{} & {} & {} & 
		$\gamma_2 \theta_{1,2} < \ell_{2,1} < \ell_{2,1} + \delta_{2,1}$ & $\gamma_2 < \ell_{2,1}$ \\
		\midrule
		% Region 7
		\multirow{2}{*}{$7$} &\multirow{2}{*}{Monostable} & \multirow{2}{*}{$\setof*{(0, 0)}$} & 
		$\ell_{1,2} < \ell_{1,2} + \delta_{1,2} < \gamma_1 \theta_{2,1}$ & 
		$\ell_{1,2} + \delta_{1,2} < 1$ \\
		{} & {} & {} & 
		$\ell_{2,1} < \ell_{2,1} + \delta_{2,1} < \gamma_2 \theta_{1,2}$ &
		$\ell_{2,1} + \delta_{2,1} < \gamma_2$ \\
		\midrule
		% Region 8
		\multirow{2}{*}{$8$} &\multirow{2}{*}{Monostable} & \multirow{2}{*}{$\setof*{(0, 1)}$} & 
		$\ell_{1,2} < \ell_{1,2} + \delta_{1,2} < \gamma_1 \theta_{2,1}$ & 
		$\ell_{1,2} + \delta_{1,2} < 1$ \\
		{} & {} & {} & 
		$\ell_{2,1} < \gamma_2 \theta_{1,2} < \ell_{2,1} + \delta_{2,1}$ &
		$\ell_{2,1} < \gamma_2 < \ell_{2,1} + \delta_{2,1}$ \\
		\midrule
		% Region 9
		\multirow{2}{*}{$9$} &\multirow{2}{*}{Monostable} & \multirow{2}{*}{$\setof*{(0, 1)}$} & 
		$\ell_{1,2} < \ell_{1,2} + \delta_{1,2} < \gamma_1 \theta_{2,1}$ & 
		$\ell_{1,2} + \delta_{1,2} < 1$ \\
		{} & {} & {} & 
		$\gamma_2 \theta_{1,2} < \ell_{2,1} < \ell_{2,1} + \delta_{2,1}$ & $\gamma_2 < \ell_{2,1}$ \\
		\bottomrule
	\end{tabular}
	\caption{The result of the DSGRN analysis for the Toggle Switch. As indicated in Figure~\ref{fig:TS}(b) DSGRN partitions the parameter space into $9$ parameter regions which give rise to the nodes of the parameter graph.  The  Morse graph associated with each node gives rise to the phenotype.
	The inequalities column indicates the explicit region in parameter space associate to each node.  
	For the reduced parameter space defined in Section \ref{sec:reducing the number}, the corresponding parameter regions are the semi-algebraic subsets of $\Xi^*$ described by the inequalities in the last column.}
	\label{tab:parameter_regions}
\end{table*}

To provide intuition with respect to parameters we note that the Hill model for the Toggle Switch is given by
\begin{equation}
    \label{eq:Toggle_ODEs}
    \begin{aligned}
    \dot{x}_1 & = -\gamma_1 x_1 + \ell_{1,2} + \delta_{1,2}\frac{\theta^{d_{1,2}}_{1,2}}{\theta^{d_{1,2}}_{1,2}+x^{d_{1,2}}_2}\\
    \dot{x}_2 & = -\gamma_2 x_2 + \ell_{2,1} + \delta_{2,1}\frac{\theta^{d_{2,1}}_{2,1}}{\theta^{d_{2,1}}_{2,1}+x^{d_{2,1}}_1}
    \end{aligned}
\end{equation}
where the parameters space is 
\begin{equation}
\label{eq:parametersHillToggle}
 \Lambda_{\text{TS}} = \setof{(\gamma_1, \ell_{1,2}, \delta_{1,2}, \theta_{1,2}, d_{1,2}, \gamma_2, \ell_{2,1}, \delta_{2,1}, \theta_{2,1}, d_{2,1})} = (0,\infty)^{10}.   
\end{equation}
We provide a more general discussion of Hill Models in Section~\ref{sec:hill models}.
Comparing the parameters listed in \eqref{eq:parametersHillToggle} with those of Table~\ref{tab:parameter_regions} there is an obvious linear projection between the parameter spaces for the Hill model and the DSGRN model for the toggle switch; $\pi_{\text{TS}}\colon \Lambda_{\text{TS}}\to \Xi_{\text{TS}}$ obtained by ignoring the Hill exponents $d_{1,2}$ and $d_{2,1}$.
As a first pass it is useful to think of DSGRN as capturing the dynamics of \eqref{eq:Toggle_ODEs} in the singular limit as the exponents $d_{1,2}$ and $d_{2,1}$ go to infinity.

In fact, there is a bijection between the equilibria that can be identified by DSGRN and the equilibria for a Hill function model for sufficiently large exponents \cite[Theorem 3.13]{duncan:gedeon:kokubu:mischaikow:oka}.

As indicated within the synthetic biology community, the Toggle Switch is of practical interest and for this community an abstract existence result invoking the phrase ``sufficiently large exponents'' is not of direct use, as typical biochemical models involve Hill exponents that are on the order of 5 or less.
This brings us back to the EQ and SNB problems described in the introduction. 
First, can we use the DSGRN information to provide intuition as to which parameter values will lead to bistability for \eqref{eq:Toggle_ODEs}?
Second, can we use this information to identify at what parameter values saddle-node bifurcations will occur? 

We note that the Morse graphs for a combinatorial network model typically encode a richer description of the global dynamical structure than simply the number and type of combinatorial attractors in each parameter region, such as the phase space region in which they occur. Analysis of the relative location of these attractors in phase space across multiple parameter regions may suggest subtle, yet global information about how these attractors continue and which types of bifurcations occur. We demonstrate this in a specific example in Section~\ref{sec:degenerate saddle-node bifurcations}.

As is done in this paper DSGRN can be used as a black box. 
However, to provide at least some intuition (see \cite{cummins:gedeon:harker:mischaikow:mok, gedeon:harker:kokubu:mischaikow:oka, duncan:gedeon:kokubu:mischaikow:oka} for details) we provide a minimal amount of information and motivation about the internal workings of DSGRN.
Returning to \eqref{eq:Toggle_ODEs} observe that as $d_{j,i}\to \infty$ the nonlinearities converge to step functions of the form
\[
\begin{cases}
\ell  + \delta &\text{if $x <\theta$} \\
\ell &\text{if $x >\theta$.}
\end{cases}
\]
Thus, in this singular limit the parameter $\theta$ indicate the locations in phase space where the dynamics changes.
DSGRN exploits this by decomposing phase space $(0,\infty)^N$ into a cubical complex defined by the hyperplanes $x_i = \theta_{j,i}$ for any edge $i$ to $j$ in the regulatory network.
DSGRN models the dynamics via a directed graph, called the \emph{state transition graph} that identifies how cells are mapped to adjacent cells (a self loop indicates that the cell is mapped to itself).
There is a single state transition graph associated to each node $s$ in the parameter graph.
DSGRN compresses the information in the state transition graph by identifying the recurrent strongly connected path components, i.e. for the ODE that is being investigated using DSGRN the strongly connected path components that contain at least one edge.
There is a 1-1 correspondence between the recurrent strongly connected path components and the elements of $\sM(s)$.
Furthermore, given $p,q\in \sM(s)$, $p\leq_s q$ implies that there exists a path in the state transition graph from the strongly connect path component associated with $q$ to the strongly connect path component associated with $p$.

Observe that if a cell in the state transition graph has a self edge, then that cell must belong to a recurrent strongly connected path component.
Furthermore, at least conceptually, the existence of a self edge suggests the existence of a fixed point. 
With this in mind  if $p\in \sM(s)$ arises from a recurrent strongly connected path component that contains exactly one cell from the cubical complex, then DSGRN labels $p$  as an $FP$, which is meant to suggest fixed point.
However, DSGRN work with purely combinatorial objects and thus it is formally incorrect to call this a fixed point.
Furthermore, it is possible that $p\in \sM(s)$ is labeled as an FP, but there is no corresponding fixed point for the ODE of interest.

If $p$ is a minimal element of $\sM(s)$ and labeled as an $FP$, then the expectation is that $p$ identifies a stable fixed point for the ODE that is being investigated using DSGRN.
Note that in this case, for any ODE for which the combinatorial DSGRN computations provide a reasonable represention there will be an associated fixed point.

As is indicated in the introduction we are interested in identifying bifurcations to differential equations. 
Since a formal definition of bifurcation involves a change in the structure of  invariant sets, DSGRN is incapable of directly identifying bifurcations.
Nevertheless changes in the Morse graphs as a function of parameter nodes strongly suggest that the dynamics for associated differential equations should also changes.
Thus, for the purpose of this paper, we introduce the following definition.

\begin{definition}\label{def:combinatorial_saddle}
Let $(\cG, p)$ be an annotated graph with a vector of interaction functions and let $\sPG(\cG, p)$ be the associated parameter graph.
A DSGRN \emph{saddle-node bifurcation edge} is an edge $(v_1,v_2)$  in $\sPG(\cG, p)$ such that, with respect to the partial order, the number of minimal FP in $\sM(v_1)$ and $\sM(v_2)$ differs by exactly one and the number of nonminimal  FP in $\sM(v_1)$ and $\sM(v_2)$ differs by exactly one.
\end{definition}

One of the aims of this paper is to develop a general procedure by which we can begin with the information from the DSGRN computations on the Toggle Switch and build up towards more complicated systems of ODEs. 
Figure~\ref{fig:strategy} is meant to provide intuition concerning the approach we take.
Consider an eight-dimensional slice of DSGRN parameter space $\Xi_{\text{TS}}$ through regions $R(4)$, $R(5)$, and $R(6)$.
As indicated in Figure~\ref{fig:TS} DSGRN identifies bistability for  $R(5)$ and monostability for $R(4)$ and $R(6)$.
For simplicity we assume $d=d_{1,2}=d_{2,1}$.
By \cite{duncan:gedeon:kokubu:mischaikow:oka} and as indicated by the gray regions the equilibria identified by DSGRN persist for $d$ sufficiently large.
For $d=1$ equation \eqref{eq:Toggle_ODEs} has a unique equilibrium.
As a consequence the bistability exhibited in $R(5)$ for large $d$ must fail for some value of $d>1$.
The generic expectation is that  bistability is lost  via a saddle-node bifurcation and thus there should be an nine-dimensional surface in $\Lambda_{\text{TS}}$ the parameter space for \eqref{eq:Toggle_ODEs} on which this bifurcation occurs.
Two simplified scenarios for the structure of this surface are shown in Figure~\ref{fig:strategy}: the solid curve and the dashed curve.
The solid curve is meant to suggest that if a saddle-node bifurcation occurs then the $\Xi_{\text{TS}}$ parameters lie in $R(5)$, whereas the dashed curve suggest that the $\Xi_{\text{TS}}$ parameters may lie in $R(4)$ or $R(6)$.
Another way of stating this is that we want to measure how saddle-node bifurcations in an explicit ODE system differ from the predicted DSGRN \emph{saddle-node bifurcation edge}.

\begin{figure}[tbh]
\begin{picture}(400,150)

\put(0,0){\begin{tikzpicture}
		
	\draw (0,0) -- (9,0);
	\draw (0,0) -- (0,4);
	\node at (9.25,0) {$\Xi$};
	
	\draw (-0.1,4) -- (0.1,4);
	\node at (-0.7,4) {$d=\infty$};
	
	\draw (-0.1,3.25) -- (0.1,3.25);
	\node at (-0.7,3.25) {$d=100$};
	\draw[dotted] (0,3.25) -- (9,3.25);
	
	\draw (-0.1,0.5) -- (0.1,0.5);
	\node at (-0.7,0.5) {$d=1$};
	
	\draw[blue] (1.75,3.25) -- (1.75,0.5);
	\node at (1.75,0.25) {${\bf r}_a$};
	\draw[blue] (3.5,3.25) -- (3.5,0.5);
	\node at (3.5,0.25) {${\bf r}_b$};
	\draw[blue] (5.5,3.25) -- (5.5,0.5);
	\node at (5.5,0.25) {${\bf r}_c$};
	
    \draw[fill=gray] (0,4) -- (3,4) -- (2.6,3.75) -- (0,3.75) -- (0,4);
    \draw[fill=gray] (3,4) -- (6,4) -- (5.9,3.75) -- (3.12,3.75) -- (3,4);
    \draw[fill=gray] (6,4) -- (9,4) -- (9,3.75) -- (6.32,3.75) -- (6,4);
    
    \draw (3,4) .. controls (4.5,0) and (4.5,0) .. (6,4);
    \draw[dashed] (3,4) .. controls (-3,0) and (12,0) .. (6,4);
    \draw[dash dot] (5.75,1.75) ellipse (0.5 and 0.2);
    
    \node at (4.5,2.5) {bistability};
%    \node at (1.25,0.75) {monostability};
    \node at (7.75,0.75) {monostability};

	\node at (1.5,-0.25) {$R(4)$};
	\node at (4.5,-0.25) {$R(5)$};
	\node at (7.5,-0.25) {$R(6)$};
	\draw (3,-0.1) -- (3,0.1);
	\draw (6,-0.1) -- (6,0.1);	
%	\draw (-0.1,4) -- (0.1,4);
			\end{tikzpicture}		
			}	
\end{picture}
\caption{Schematic description of the strategy of this paper. $R(4)$, $R(5)$, and $R(6)$ are regions of parameter space $\Xi$ identified by DSGRN for Toggle Switch.
As indicated in Figure~\ref{fig:TS} DSGRN identifies bistability for  $R(5)$ and monostability for $R(4)$ and $R(6)$.
We are interested in understanding the equilibrium structure for \eqref{eq:Toggle_ODEs}.
Set $d=d_{1,2}=d_{2,1}$.
By \cite{duncan:gedeon:kokubu:mischaikow:oka} and as indicated by the gray regions the equilibria identified by DSGRN persist for $d$ sufficiently large.
For $d=1$ \eqref{eq:Toggle_ODEs} has a unique equilibrium.
There therefore must be an interface between the monostability region and the bistability region. In this paper, we use paths of the form of ${\bf r}$ to justify that the main mechanism to pass from monostability to bistability is through a saddle node bifurcation and we investigate the ``shape'' of the boundary between stability regions. 
}
\label{fig:strategy}
\end{figure}

To test which of these scenarios is more appropriate we consider curves through parameter space parameterized by $d$.
Given a curve of the form $\curve_a$ we expect to see no saddle-node bifurcation if the solid curve applies and two saddle-node bifurcations if the dashed curve applies. 
In addition, if the solid curve applies then we expect that the value of $d$ at which the saddle-node bifurcation occurs
for curves of the form $\curve_b$ will decrease as a function of the distance the $X_i$ parameter is from the boundary of $R(4)$ or $R(5)$.
As indicated in Section~\ref{sec:toggle_switch} it appears that the solid curve provides the better intuition. What is perhaps more surprising is that there exist curves $\curve_c$ along which isolas (demarcated by the dash-dot curve) of equilibria appear, i.e., as $d$ is decreased two equilibrium appear and then the same two equilibrium disappear via a saddle-node bifurcation.

\section{Hill models}
\label{sec:hill models}
%!TEX root = dynamics_main.tex

In this section we provide a formal definition of the Hill models 
considered in this paper.
We begin by noting that Hill models take the form 
\begin{equation}
	\label{eq:Hill_model}
	\dot x = -\Gamma x + \cH(x), \qquad x\in [0,\infty)^N
\end{equation}
where $\Gamma$ is a diagonal matrix with $\Gamma_{ii} = \gamma_i >0$.
The rest of this section is dedicated to defining the nonlinear function $\cH$.
This requires several preliminary constructions that are imposed so that there is an obvious correspondence between the computations that DSGRN is currently capable of performing  and the form of the Hill models. 

\begin{definition}
\label{def:hill_function_response}
An \emph{activating Hill response} is a function $H^+ : [0, \infty) \to [0, \infty)$ given by 
\[
H^+\paren*{x} := \ell + \delta \frac{x^\hill}{\theta^\hill + x^\hill} 
\]
and a \emph{repressing Hill response} is a function $H^- : [0, \infty) \to [0, \infty)$ given by
\[
H^{-}\paren*{x} := \ell + \delta \frac{\theta^\hill}{\theta^\hill + x^\hill} 
\]
where the parameters $\setof{\ell, \delta, \theta, \hill}$ are positive  and $\hill \geq 1$. We refer to the parameter $\hill$ as the \emph{Hill exponent}.
\end{definition}

In biological applications, Hill models are extensively used for modeling the internal behaviour of the cell~\cite{collins,levasseur1998modeling,yan2001theoretical}. They are coupled with the study of network dynamics, as presented in~\cite{DSGRN}, where a repressing edge is associated with $H^-$ and an activating edge with $H^-$. Observe that $H^+$ is monotonically increasing and satisfies 
\[
H^+(0) = \ell \qquad \lim\limits_{x \to \infty} H^+(x) = \ell + \delta.
\]
Similarly, $H^-$ is monotonically decreasing and satisfies 
\[
H^{-}(0) = \ell + \delta \qquad  \lim\limits_{x \to \infty} H^-(x) = \ell.
\]

In general $\cH$ is a function of multiple variables and in what follows we indicate how $\cH$ can be constructed from activating and repressing Hill responses.
%>>

\begin{definition}
	\label{def:interaction_function}
	A polynomial, $p \in \rr[z_1, \dotsc, z_N]$ is called an {\em interaction function} if it has the form
	\[
	p = \prod_{m = 1}^{q} p_m, \qquad \text{where} \quad  p_m = \sum_{j \in I_m} z_j,
	\]
	and $\setof{I_1, \dotsc, I_q}$ is a partition of the integers, $\setof{1, \dotsc, N}$. 
\end{definition}

Interaction functions allow us to define Hill models.

\begin{definition}
\label{def:HillModel}
A \emph{Hill model} is a differential equation $\dot{x}=f(x)$, $x\in[0,\infty)^N$ such that for each  $n \in \setof{1, \dotsc, N}$,
\begin{equation}
	\label{eq:hill_model_coordinate}
	\dot{x}_n = f_n(x) = -\gamma_n x_n + \cH_n(x)
    \end{equation}
where
    \begin{equation}
	\label{eq:nonlinear_coordinate}
	\cH_n(x) = p_n \paren*{H^*_{n,1}(x_1), \dotsc, H^*_{n,N}(x_N)}
    \end{equation}
for some interaction function $p_i$ and $*$ indicates $+$ or $-$.
The associated parameter space is denoted by $\Lambda$.
\end{definition}

\begin{remark}
\label{rem:projection}
To match a Hill model with DSGRN, one interprets the edges of the annotated graph $(\cG,p)$ as
    \[
    H^*_{i,j} = \begin{cases}
    H^+_{i,j} & \text{if $j\to i$}\\
    H^-_{i,j} & \text{if $j\dashv i$}\\
    0 & \text{otherwise.}
    \end{cases}
    \]
The Hill models described in Definition~\ref{def:HillModel} are based on the DSGRN parameters $\Lambda$ plus the Hill exponents $\hill\in(0,\infty)^E$. 
We let $\pi\colon \Lambda \to \Xi$ denote the linear projection between these two parameter spaces obtained by ignoring the Hill exponents.
\end{remark}

The following proposition suggests that the global  dynamics of Hill systems are trivial for sufficiently small Hill exponents. 
\begin{proposition}
\label{prop:d=0}
Consider a Hill model where each Hill exponent satisfies $d_{n,m}=0$.
Then, there exists a unique equilibrium that is a global attractor.
\end{proposition}
\begin{proof}
Observe that if $d_{n,m}=0$, then
\[
H_n^\pm(x_m) = \ell_{n,m} + \frac{1}{2}\delta_{n,m} 
\]
and hence is a positive constant.
Thus, \eqref{eq:Hill_model} takes the form
\[
\dot{x} = -\Gamma x + P
\]
where $P$ is a constant vector with positive entries.
\end{proof}

Though suppressed in the notation presented above, Hill models involve a multitude of parameters.
In particular, each state variable $x_i$ contributes a linear decay parameter $\gamma_i$ and if  $x_j$ influences the rate of change of $x_i$, then the associated Hill function contributes four parameters $\ell_{i,j}$, $\delta_{i,j}$, $\theta_{i,j}$, and $\hill_{i,j}$. For fixed parameter $\lambda$ we will write the Hill model as
\[
\dot{x} = f(x,\lambda).
\]

\section{Equilibria and saddle-node bifurcations}
\label{sec:equilibria}
%!TEX root = dynamics_main.tex
Our strategy for identifying saddle-node bifurcations for \eqref{eq:abstract} involves four steps.
\begin{description}
    \item[Step 1] Choose a parameterized path $\curve\colon [0,1]\to \Lambda$ along which we expect a saddle-node bifurcation to occur. 
    \item[Step 2] Subdivide $[0,1]$ at points $0=s_0 <\cdots < s_i <\cdots < s_I = 1$ and identify equilibria of \eqref{eq:abstract} at parameter value $\curve(s_i)$.
    
    \item[Step 3] If the number of equilibria for $f$ at $\curve(s_i)$ and $\curve(s_{i+1})$ differs, 
    then a refined search takes place to precisely determine where a bifurcation should take place.

    \item[Step 4] 
    Using candidates identified in {\bf Step 3}  numerically identify saddle-node bifurcations.
\end{description}
We begin with a few high level comments concerning these four steps.

We remark that the numerical capabilities of NDMA can be applied to any path $\curve\colon [0,1]\to \Lambda$. 
However, for the sake of simplicity throughout the remainder of this paper we restrict our attention to the following two types of paths. 

Define $\Lambda_\bullet \subset \Lambda$ by equating all Hill exponents, i.e., $d_{i,j}= d_{n,m}$ for all $i,j,n,m = 1,\ldots, N$.
With regard to {\bf Step 1} we choose $\curve\colon [0,1]\to \Lambda_\bullet$ such that $\pi(\curve(s))\in \Xi$ is constant for all $s$. 
Furthermore, unless otherwise stated we assume that the Hill exponent for $\curve (s)$ is given by $d=1 + 99s$, i.e., the value of the Hill exponent increases linearly as a function of $s$ starting at $d=1$ and ending at $d=100$. We call this a \emph{vertical path}. 

For the second type of path, the Hill exponents remain unchanged along the path and we refer to it as a \emph{horizontal path}. The justification for these two types of paths is in the combinatorial background just introduced. Indeed, a vertical path will always remain in the given combinatorial parameter region, while a horizontal path could cross the boundary between parameter regions previously defined. 
%>>

Identification of equilibria ({\bf Step 2}) is a well studied but nontrivial problem and is the focus of most of this section.
There are two numerical questions that we need to address.
The first, where are the equilibria? 
We answer this in general in Section~\ref{sec:overviewEquilibrium} and more specifically in Section~\ref{sec:the bootstrap algorithm} where we exploit the structure of Hill models and provide an algorithm (see Algorithm~\ref{alg:bootstrap_equilibria}) that identifies regions in phase space in which the equilibria must occur.
The second, how to determine whether different numerical representations identify distinct equilibria?
We address this in Section~\ref{sec:isolating equilibria}.

For typical applications, one might expect that {\bf Step 3} involves the use of a continuation algorithm for equilibria.
In particular, a continuation algorithm based on pseudo-arc length would simultaneously detect saddle-node bifurcations.
However, our choice of parameter path, based on the Hill exponent, leads to a stiff numerical problem: far from a saddle node bifurcation, the equilibria move little in phase space, while close to saddle nodes significant changes occur depending on small changes of the Hill exponent.
For this reason we employ a bisection algorithm that appears to be not only more reliable, but also faster than its continuation counterpart.

Turning to {\bf Step 4}, we refer the reader to \cite[Section 8.2]{chicone}  and \cite[Theorem 8.12]{chicone} for a definition of a saddle-node bifurcation and analytic conditions that guarantee its existence, respectively.
The following theorem (see \cite{lessard}) provides a tool by which we can numerically identify that a saddle-node bifurcation occurs.

\begin{theorem}
\label{thm:saddle_node_bifurcation}
Consider a smooth one parameter family of ODEs, 
\[
\dot x = g(x, s), \quad x\in \rr^N, s\in\rr .
\]
Define $G : \rr^{2N+1} \to \rr^{2N+1}$ by
\begin{equation}\label{eq:num_saddle_node}
G(x, v, s) := 
\begin{pmatrix}
	g(x, s)\\
	D_xg(x, s) v\\
	v^Tv - 1 
\end{pmatrix}
\qquad x, v \in \rr^N, s \in \rr.
\end{equation}
Let $u = (x, v, s) \in \rr^{2N+1}$ and suppose $\hat{u} := (\hat{x}, \hat{v}, \hat{s})$ is a root of $G$ satisfying
\begin{enumerate}
	\item $D_u G (\hat u)$ is an isomorphism. 
	\item Every nonzero eigenvalue of $D_x g(\hat x, \hat{s})$ has nonzero real part. \end{enumerate}
Then $g$ undergoes a saddle-node bifurcation at $(\hat{x}, \hat{s})$ and $\ker D_x g(\hat{x}, \hat{s}) = \inspan{\setof{\hat{v}}}$.
\end{theorem}

\begin{remark}
We emphasize that the above mentioned steps are implemented in a pure numerical fashion in this paper.
Application of Theorem~\ref{thm:saddle_node_bifurcation} in conjunction with
validation techniques would result in {\emph{proofs}} of the existence of a saddle node bifurcation. 
We do not in this paper pursue these ideas further, but the interested reader can refer to \cite{lessard} for additional details. 
We also note that all the computational elements needed for the validation of a Hill model are included in \library.
\end{remark}

\subsection{Overview of equilibrium finding}
\label{sec:overviewEquilibrium}
Within {\bf Step 2} we are asked to numerically identify  equilibria for a given parameter value. Conceptually our approach involves two black box algorithms: {\tt FindRoot} and {\tt Unique}.
The first, takes a function, $f : \rr^N \to \rr^N$, and an initial guess $x_0 \in \rr^N$ as input and when successful returns $\hat x$ satisfying $\norm{f(\hat x)} \approx 0$. 

The second takes pairs of distinct outputs from {\tt FindRoot} and determines heuristically whether or not they represent different approximations for the same root of $f$ or distinct roots.
More specifically, assume there exists $\bar{x}$ and $\epsilon >0$ such that  $f(\bar{x})=0$ and $0<\| x-\bar{x}\|<\epsilon$ implies $f(x)\neq 0$.
If $\hat x_1, \hat x_2$ satisfy 
\[
\norm{f(\hat{x}_1)} \approx 0 \approx \norm{f(\hat x_2)} \qquad \text{and} \qquad \norm{\bar{x} - \hat{x}_i} < \epsilon,\ i=1,2
\] 
then {\tt Unique}$(f,\hat x_1, \hat x_2) = \hat x_1$. 
Otherwise, {\tt Unique}$(f,\hat x_1, \hat x_2) = \setof{\hat x_1,\hat{x}_2}$ and we say that $\hat{x}_1$ and $\hat{x}_2$ are \emph{approximately} distinct with respect to $f$. In Section \ref{sec:isolating equilibria} we provide details about how the parameter $\epsilon$ can be reliably chosen. 

More generally, if $\hat x$ is an array of vectors in $\rr^N$ and $f$ a function from $\rr^N$ to itself, then {\tt Unique}$(f,\hat x)$ returns an new array of vectors in $\rr^N$ of size equal or smaller than $\hat x$, in which each pair of vectors is approximately distinct.

It is worth noting that Hill functions are trivially bounded, thus for any Hill model, there exists a  rectangular subset of $X$ of the form 
\begin{equation}
    \label{eq:Rectangle}
R := \prod_{i = 1}^{N} [a_i, b_i], \qquad [a_i, b_i] \subset (0, \infty) \forall 1 \leq i \leq N.
\end{equation}
such that all the zeros of the Hill model are within this rectangle. Furthermore, computing the bounds of $R$ is trivial. 
The details are presented in Section \ref{sec:the bootstrap algorithm}.

Algorithm \ref{alg:general_equilibria} indicates how given $R$ defined by \eqref{eq:Rectangle} we utilize {\tt FindRoot} and  {\tt Unique} to search for zeros of the Hill model $f$ at a fixed parameter value.
In words, each interval in the product \eqref{eq:Rectangle} is partitioned into $k$ subintervals bounded by $k+1$ uniformly spaced nodes. 
The product of these nodes forms a grid of points in $\rr^N$ which covers $R$. 
Each of the $(k+1)^N$ points in this grid is taken as an initial condition for {\tt FindRoot} which attempts to return a candidate equilibrium nearby. 
The algorithm returns an array containing such candidates which are not identified as equivalent by {\tt Unique}.
Increasing $k$ increases the likelihood of finding all equilibria.

\begin{algorithm}
	\caption{General algorithm}
	\label{alg:general_equilibria}
	\begin{algorithmic}[1] % The number tells where the line numbering should start
		\Function{{\tt HillEquilibria}}{$f, R, k$} 
		\State $\hat{x} \gets ()$ \Comment{Initialize equilibrium array}
		\State $\Delta_i \gets \frac{b_i - a_i}{k}$
		\State $u_i \gets (a_i, a_i + \Delta_i, \dotsc, a_i + (k-1) \Delta i, b_i)$ \Comment{Discretize factors}
		\For{$\kappa \in \setof{1,\dotsc, k}^N$} 
		\State $x_0 \gets (u_{1, \kappa_1}, \dotsc, u_{N, \kappa_N})$ 
		\State $r_\kappa \gets {\tt FindRoot}(f, x_0)$ \Comment {Returns a candidate when it converges}
		\If{$r_\kappa$}
		%$\hat x(:, -1) := r_\kappa$
		$\hat x. {\tt Append}(r_k) $
		 \Comment {Append candidate equilibrium}
		\EndIf
		\EndFor
		\State \textbf{return} {\tt Unique}$(f, \hat x) $
		\EndFunction
	\end{algorithmic}
\end{algorithm}

\subsection{The boxy box algorithm}
\label{sec:the bootstrap algorithm}

In this section we define an algorithm that exploits the structure of Hill models in specific cases to efficiently localize equilibria.
The main idea is to begin with an initial rectangular subset  of $(0,\infty)^N$ which is an enclosure for all equilibria and then iteratively obtain tighter rectangular enclosures.

\begin{definition}
	\label{def:monotone_factorization}
A continuous function $g\colon [0,\infty)^N\to (0,\infty)^N$ has a {\em monotone factorization} if for each $i=1,\ldots, N$, $g_i$ factors into a product of the form
\[
g_i(x) = g_i^+(x) g_i^-(x) \quad \forall x \in [0, \infty)^N,
\]
where $g_i^+\colon [0,\infty)^N\to (0,\infty)^N$ is bounded and strictly increasing with respect to $x_1,\dotsc, x_N$ and similarly, $g_i^-\colon [0,\infty)^N\to (0,\infty)^N$ is bounded and strictly decreasing with respect to $x_1, \dotsc, x_N$. 
\end{definition}

Consider a continuous function $f\colon [0,\infty)^N\to\rr^N$ of the form $f(x) = -\Gamma x +g(x)$ where
$g\colon [0,\infty)^N\to(0,\infty)^N$  has a monotone factorization and for all $i=1,\ldots, N$, $f_i(x) = -\gamma_i x + g_i(x)$ and $\gamma_i >0$.
Define  $\Phi\colon \rr^{2N}\to \rr^{2N}$ coordinate-wise by the formulas
\[
		\Phi_i(\concat(\alpha, \beta)) = \frac{1}{\gamma_i} g_i^+(\alpha) g_i^-({\beta}) \quad\text{and}\quad
		\Phi_{N+i}(\concat(\alpha, \beta)) = \frac{1}{\gamma_i} g_i^+(\beta) g_i^-(\alpha), \qquad i = 1, \ldots, N
\]
where $\alpha, \beta \in \rr^N$, and where $\concat(\alpha, \beta)\in \rr^{2N}$ is the standard concatenation of two vectors.
\begin{theorem}
	\label{thm:bootstrap_eqbounds}
	Consider $f$ and $\Phi$ as defined above. 	
	Assume that $\liminf\limits_{\|x\|\to \infty}g_i^-({x})>0$ for all $i=1,\ldots,N$.
	Then, the following are true. 
	\begin{enumerate}[align = Center, label=\roman*]
	\item $x\in [0,\infty)^N$ is a zero of $f$ if and only if $\concat(x,x) \in [0, \infty)^{2N}$ is a fixed point of $\Phi$. 
    \item Define $\concat\paren*{\alpha^{(0)},\beta^{(0)}} \in \rr^{2N}$ coordinate-wise by
    \begin{equation}
    \alpha^{(0)}_i := \frac{1}{\gamma_i} g_i^+(0) \liminf_{\|x\|\to \infty}g_i^-({x})
    \quad\text{and}\quad
    \beta^{(0)}_i :=  \frac{1}{\gamma_i} \limsup_{\|x\|\to \infty}g_i^+(x) g_i^-(0)
    \end{equation}
    and iteratively define $\concat(\alpha^{n+1},\beta^{n+1}) = \Phi(\concat(\alpha^{n},\beta^{n}))$ for $n \geq 1$. Then,  $\concat(\hat{\alpha},\hat{\beta}) := \lim_{n\to\infty}\concat(\alpha^{n},\beta^{n})$ exists.
    \item If $f(\hat x)=0$, then 
    \[
    \hat{\alpha}_i \leq \hat{x}_i\leq \hat{\beta}_i, \quad \forall i=1,\ldots,N.
    \]
    \end{enumerate}
\end{theorem}

\begin{proof}
    We leave it to the reader to check ({\it{i}}). ({\it{ii}}) follows from the boundedness and strict monotonicity of $g_i^+$ and $g_i^-$. To be more specific, we prove inductively that for $1 \leq i \leq N$, $\alpha_i^{(n)}$ and $\beta_i^{(n)}$ are monotonically increasing and decreasing sequences, respectively. 
    The base case follows from easy estimates
    \[
   \alpha_i^{(1)} = \Phi_i(\concat\paren*{\alpha^{(0)},\beta^{(0)}})  = \frac{1}{\gamma_i} g_i^+\paren*{\alpha^{(0)}} g_i^-\paren*{\beta^{(0)}} 
		 > \frac{1}{\gamma_i} g_i^+\paren*{0} \liminf_{\|x\|\to \infty}g_i^-({x}) 
		 = \alpha^{(0)}_i
	\]
and
	\[
	\beta_i^{(1)} =	\Phi_{N+i}(\concat\paren*{\alpha^{(0)},\beta^{(0)}})  = \frac{1}{\gamma_i} g_i^+\paren*{\beta^{(0)}} g_i^-\paren*{\alpha^{(0)}}  
		 < \frac{1}{\gamma_i} \limsup_{\|x\|\to \infty}g_i^+\paren*{x} g_i^-\paren*{0} 
		 = \beta^{(0)}_i.
	\]
    
    Now assume that $\alpha_i^{(n)} < \alpha_i^{(n-1)}$ and $\beta_i^{(n)}> \beta_i^{(n-1)}$.  
    The strict monotonicity of $g_i^+$ and $g_i^-$ implies that
    \[
   \alpha_i^{(n+1)} = \Phi_i(\concat\paren*{\alpha^{(n)},\beta^{(n)}})  = \frac{1}{\gamma_i} g_i^+\paren*{\alpha^{(n)}} g_i^-\paren*{\beta^{(n)}} 
		 > \frac{1}{\gamma_i} g_i^+\paren*{\alpha^{(n-1)}} g_i^-(\beta^{(n-1)}) 
		 = \alpha^{(n)}_i
	\]
and
	\[
	\beta_i^{(n+1)} =	\Phi_{N+i}(\concat\paren*{\alpha^{(n)},\beta^{(n)}})  = \frac{1}{\gamma_i} g_i^+\paren*{\beta^{(n)}} g_i^-\paren*{\alpha^{(n)}}  
		 < \frac{1}{\gamma_i} g_i^+\paren*{\beta^{(n-1)}} g_i^-\paren*{\alpha^{(n-1)}} 
		 = \beta^{(n)}_i.
	\]
	
	The proof of ({\it{iii}}) is also done inductively.
	Define 
	\[
    \cR^{(n)} := \prod_{i = 1}^{N} [\alpha^{(n)}_i, \beta^{(n)}_i].
    \]
    By the proof of ({\it{ii}}), $\cR^{(n+1)}\subset \cR^{(n)}$.
	Define  $F : [0,\infty)^N \to [0,\infty)^N$ by the formula
\[
F_i(x) = \frac{1}{\gamma_i }g_i(x) \qquad 1 \leq i \leq N.
\]
Observe that if $f(\hat{x})=0$, then $F(\hat x)=\hat{x}$. 
Therefore, it suffices to prove that if $F(\hat x)=\hat{x}$ then  $\hat{x}\in \cR^{(n)}$ for all $n \in \nn$. 

Observe that from the definitions of $F$ and $\paren*{\alpha^{(0))},\beta^{(0)}}$, 
\[
F\paren*{[0,\infty)^N} = \cR^{(0)}
\]
and therefore $\hat{x} \in \cR^{(0)}$. 

Inductively, suppose that $n \in \nn$ is fixed and $\hat{x}  \in \cR^{(n-1)}$, i.e.
\begin{equation}
	\label{eq:induction_bound}
\alpha_i^{(n-1)} \leq \hat{x}_i \leq \beta_i^{(n-1)} \qquad \text{for} \ 1 \leq i \leq N. 
\end{equation}
The inequalities of \eqref{eq:induction_bound} combined with the definition of $\Phi$ implies that for all $1 \leq i \leq N$
\[
\alpha_i^{(n)} = \frac{1}{\gamma_i} g_i^+\paren*{\alpha^{(n-1)}} g_i^-\paren*{\beta^{(n-1)}} \leq F_i(\hat x)=\hat{x}_i \leq \frac{1}{\gamma_i} g_i^+\paren*{\beta^{(n-1)}} g_i^-\paren*{\alpha^{(n-1)}} = \beta_i^{(n)}
\]
where the inequalities are obtain from the fact that $g_i^+$ and $g_i^-$ are strictly monotonically increasing and decreasing respectively.
Therefore, $\hat{x} \in \cR^{(n)}$.
\end{proof}

Observe that $F\paren*{\cR^{(n-1)}} = \cR^{(n)}$.

 This motivates the following algorithm for bounding equilibria. 
\begin{algorithm}
	\caption{boxy box algorithm}
	\label{alg:bootstrap_equilibria}
	\begin{algorithmic}[1] % The number tells where the line numbering should start
				\Function{\tt RootEnclosure}{$f$} 
		\State $u \gets (\alpha^{(0)}, \beta^{(0)})$ \Comment{Initialize orbit as described in Theorem \ref{thm:bootstrap_eqbounds}}
		\State $v \gets \Phi(u)$ 
		\While{$\norm{u - v} > \epsilon$} 
		\State $u \gets v$
		\State $v \gets \Phi(u)$
		\EndWhile 
		\State $R \gets \prod_{i = 1}^N [v_i, v_{N+i}]$
		\State \textbf{return} $R$
		\EndFunction
		\Function{\tt MonotoneHillEquilibria}{$f, k$}
		\State $R \gets ${\tt RootEnclosure($f$)}
		\State \textbf{return} {\tt HillEquilibria$(f, R, k)$}
		\EndFunction
	\end{algorithmic}
\end{algorithm}

\begin{remark}\label{rem:saddlenode_alg}
The boxy box algorithm  can be used to detect possible bifurcations when applied to a parameter dependent problem along a path of parameters. Indeed, let $g_p$ be a parameter dependent continuous function with a monotone factorisation. Furthermore, let $p_1$ and $p_2$ be two parameters, such that $\|p_1 - p_2\|<< 1$. Then, if the boxy box algorithm returns a degenerate box $\tilde x$ for $g_{p_1}$ and a non-degenerate box $[x_1, x_2]$ for $g_{p_2}$, then we want to further investigate the behaviour since such change in dynamics is often indicative of a bifurcation.

While studying saddle node bifurcations,  we noticed that one of the corners of the newly formed box corresponds to the previously found equilibrium, while the opposite corner tightly encloses a new equilibrium. 
In our numerical experiments, this phenomenon was observed in the overwhelming majority of cases.
\end{remark}

\begin{remark}
    In is important to underline that the boxy box algorithm might (and sometimes does) converge to a box that is not a tight enclosure of the equilibria of the system, these cases are exceedingly rare in the situations at hand. Thanks to this numerical stability, applying the boxy box algorithm to a Hill system allows us to immediately detect both monostability and multistability. Furthermore, it is often the case that equilibria are close to the corners of the box (i.e., the equilibrium has coordinates either $\alpha_i$ or $\beta_i$ for every component $i$), that allows for fast and stable convergence of any other equilibrium search algorithm initialised on only the corners of the box. Considering these two elements, the application of this algorithm slashes the computational time for equilibrium searches in all our test cases, such that plots such as Figure \ref{fig:coherency_TS} take mere minutes, while more challenging ones like Figure \ref{fig:coherencyVSneqsANDhill} take around 2 hours on a laptop (MacBook Pro 2019 2.3 GHz1117 8-Core Intel Core i9).
\end{remark}

\subsection{Isolating equilibria}
\label{sec:isolating equilibria}

In this section we provide the algorithm {\tt Unique} required for Algorithm~\ref{alg:general_equilibria}.
The key theoretical result is the following theorem \cite[Chapter~1]{VandenBerg2018}.

\begin{theorem}
\label{thm:rpolynomial}
Let $F\colon \rr^N\to \rr^N$ be a differentiable function and
assume that $\hat{x}\in \rr^N$ is a numerical zero of $F$, i.e., $F(\hat{x})\approx 0$.
Let $A$ be an $n\times n$ matrix.
Suppose
\[
Y \geq \|AF(\hat{0}) \|,\quad Z_0 \geq \| I-ADF(\hat{x})\|, \quad\text{and}\quad Z_1(r) \geq \max_{z,b\in B_1(0)}\|AD^2F(\hat x + rz)b\|.
\]
Define the radii polynomial by
\begin{equation}\label{eq:radii_pol}
p(r) = Z_1(r)r^2 +(Z_0  -1)r +Y.
\end{equation}
If there exists $r^*$ such that $p(r^*)<0$, then there exists a unique $\bar{x}$ such that $F(\bar{x}) =0$ and $\|\hat{x} -\bar{x}\| < r^*$.
\end{theorem}

On the assumption that the second derivative of $F$ is essentially constant on a small ball we approximate $Z_1(r)$ by the constant $Z_1 :=  \|AD^2F(\hat x)\|$.
This reduces the radii polynomial to 
\[
p(r) = Z_1r^2 +(Z_0  -1)r +Y
\]
for which the biggest existence radius is 
\[
r_{\hat{x}}^* = \frac{1-Z_0+ \sqrt{Z_1^2 - 4(Z_0-1) Y}}{2Z_1}.
\]

We now return to the question of identifying distinct zeros, Given numerical zeros, $\hat{x}_1$ and $\hat{x}_2$ of $F$, compute  $r_{\hat{x}_1}^*$ and $r_{\hat{x}_2}^*$.
If
\[
\|\hat{x}_1 - \hat{x}_2 \| < \max \{ r_{\hat{x}_1}, r_{\hat{x}_2}\}
\]
then we consider them to be numerical representations of the same zero, otherwise they represent different zeros. 

More generally, given a set of approximate zeros $\{\hat x_1, \dots \hat x_n\}$ compute $r^*_{\hat{x}_i}, i=1,\dots, n$, and then delete any $x_i$ such that there is an $x_j$ with $\| x_i - x_j\| < \max( r^*_{\hat{x}_i}, r^*_{\hat{x}_j})$. 
To reduce the number of comparisons, we restrict ourselves to $j<i$.

We implement \texttt{RadiiPol} as a function that computes $r^*$ for any function $f$ and an associated approximate equilibrium $\hat x$.

\begin{algorithm}
	\caption{Unique algorithm}
	\label{alg:unique_equilibria}
	\begin{algorithmic}[1] % The number tells where the line numbering should start
				\Function{\tt Unique}{$f$, $X$}
		\For{$\hat x\in X$}
		    \State $r^* \gets ${\tt RadiiPol}$(f, \hat x)$
		    \If{ $\hat y\in X$: $\|\hat x-\hat y\|<r^*$}
		        \State delete $\hat y$ from $X$
		    \EndIf
		\EndFor
		\State \textbf{return} $X$
		\EndFunction
	\end{algorithmic}
\end{algorithm}

\section{Bifurcation analysis of the Toggle Switch}
\label{sec:toggle_switch}
%!TEX root = dynamics_main.tex

In this section we discuss the dynamical behaviour of the Toggle Switch depending on parameters, studying both the presence of equilibria and their appearance through bifurcations - namely saddle node bifurcations. 
At first, we present some techniques to reduce the number of parameters used in the Toggle Switch, then we used such reduction to make the results more graphically compelling. With these rewritings, we can approach the numerical results more directly. We can then discuss which high-codimension bifurcations could be hidden behind the results we are experiencing.

\subsection{Parameter reduction}
\label{sec:reducing the number}

While  the full parameter space of the Toggle Switch can be studied by \library, we are interested in making the results more readable and easier to plot. For this, we will make several changes to the model in order to reduce the dimension of the parameter space via non-dimensionalization of the parameters. 
We may fix $3$ parameter values to be $\gamma_1 = \theta_{2,1} = \theta_{1,2} = 1$, then we obtain the {\em reduced Toggle Switch} Hill model defined by 
\begin{equation}
    f^*(x) :=   
\begin{pmatrix}
	- x_1 + \ell_{1,2}^* + \delta_{1,2}^* \frac{1}{1 + x_2^\hill} \\
	- \gamma_2^* x_2 + \ell_{2,1}^* + \delta_{2,1}^* \frac{1}{1 + x_1^\hill} 
\end{pmatrix}.
\end{equation}
where the parameters of the reduced model are related to the original parameters by the identities
\begin{align*}
	\ell_{1,2}^* & := \frac{\ell_{1,2}}{\gamma_1 \theta_{2,1}}, \qquad
	\delta_{1,2}^* & := \frac{\delta_{1,2}}{\gamma_1 \theta_{2,1}}, \qquad
	\ell_{2,1}^* & := \frac{\ell_{2,1}}{\gamma_1 \theta_{1,2}}, \qquad
	\delta_{2,1}^* & := \frac{\delta_{2,1}}{\gamma_1 \theta_{1,2}}, \qquad
	\gamma_2^* & := \frac{\gamma_2}{\gamma_1}.
\end{align*}
The parameter space associated to $f^*$ is the subspace, $\Lambda^*_\bullet \subset \Lambda_\bullet$, defined by
\[
\Lambda^*_\bullet := \setof*{\lambda \in \Lambda : \gamma_1 = \theta_{2,1} = \theta_{1,2} = 1, \ \hill_{1,2} = \hill_{2,1} = d} \cong (0, \infty)^6,
\]
and we denote a typical parameter by
\[
\lambda^* := \paren*{\ell_{1,2}^*, \delta_{1,2}^*, \gamma_2^*, \ell_{2,1}^*, \delta_{2,1}^*, d} \in \Lambda^*_\bullet.
\] 

All computations described in the remaining sections are performed using the reduced Toggle Switch model. 
This allows us to take advantage of the reduced number of parameters for faster computation.
We note that the {\library} library has been written to allow these sorts of constraints to be implemented just as easily as a general Hill model.

After imposing these parameter constraints for the DSGRN parameter regions we obtain a {\em reduced} version of the combinatorial parameter space denoted by $\Xi^* := (0, \infty)^5 \subset \rr^5$ where a typical reduced combinatorial parameter has the form $\xi^* = \paren*{\ell_{1,2}^*, \delta_{1,2}^*, \gamma_2^*, \ell_{2,1}^*, \delta_{2,1}^*}$. Observe that the $9$ DSGRN parameter regions are projected onto semi-algebraic subsets of $\Xi^*$ defined by the {\em reduced inequalities} given in the last column of Table \ref{tab:parameter_regions}.  Analogous to the discussion in Remark~\ref{rem:projection} we define a projection map for the reduced parameter space $\pi_{\Xi^*} : \Lambda^*_\bullet \to \Xi^*$ which is defined by projection onto the first $5$ coordinates.

\subsection{The boxier box algorithm for the Toggle Switch} \label{sec:numerical_analysisTS}

The boxy box algorithm gives us additional information on the dynamics when applied to the Toggle Switch. In this section we present these results, summarised in Theorem \ref{thm:toggle_bootstrap_eqbounds}. 

We begin by observing that $f^*$ satisfies Definition \ref{def:monotone_factorization} for any $\lambda^* \in \Lambda^*$ and therefore the boxy box algorithm is applicable. Following the construction in Section \ref{sec:the bootstrap algorithm} we obtain the bootstrap map for the reduced Toggle Switch which is given by the formula
\[
\Phi\left(\concat(\alpha, \beta)\right) = 
\begin{pmatrix}
	\ell_{1,2}^*+ \frac{\delta_{1,2}^*}{1 + \beta_2^\hill} \\
	\frac{1}{\gamma_2^*} \paren*{\ell_{2,1}^*+ \frac{\delta_{2,1}^*}{1 + \beta_1^\hill}} \\
	\ell_{1,2}^*+ \frac{\delta_{1,2}^*}{1 + \alpha_2^\hill}\\
	\frac{1}{\gamma_2^*} \paren*{\ell_{2,1}^* + \frac{ \delta_{2,1}^*}{1 + \alpha_1^\hill}}
\end{pmatrix}
\qquad \begin{pmatrix}\alpha,\beta \end{pmatrix}\in \rr^4. 
\]
Following Algorithm \ref{alg:bootstrap_equilibria} we start with the initial condition 
\[
\concat\begin{pmatrix}\alpha^{(0)}, \beta^{(0)} \end{pmatrix}:= 
\begin{pmatrix}
\ell_{1,2}^*\\
\ell_{2,1}^* \\
\ell_{1,2}^*+ \delta_{1,2}^*\\
\ell_{2,1}^* + \delta_{2,1}^*
\end{pmatrix},
\]
and define the iterates under $\Phi$ by $\begin{pmatrix}\alpha^{(n)}, \beta^{(n)} \end{pmatrix} = \Phi\left(\concat\begin{pmatrix}\alpha^{(n-1)}, \beta^{(n-1)}\end{pmatrix}\right)$ for $n \geq 1$. Theorem \ref{thm:bootstrap_eqbounds} ensures that the sequence of iterates converges to a fixed point of $\Phi$ denoted by $
\concat\begin{pmatrix}\hat\alpha, \hat\beta \end{pmatrix} \in \rr^4$ and that all equilibria of $f^*$ are contained in the rectangle $[\hat{\alpha}_1, \hat{\beta}_1] \times [\hat{\alpha}_2, \hat{\beta}_2] \subset X$. However, in the case of the Toggle Switch we can prove a stronger version of Theorem \ref{thm:bootstrap_eqbounds} which we exploit in our numerical implementation.

\begin{theorem}
	\label{thm:toggle_bootstrap_eqbounds}
	Let $f^*$ denote the Toggle Switch Hill model with $\lambda^* \in \Lambda^*$ fixed and let $\Phi \colon \rr^4 \to \rr^4$ be the associated bootstrap map for $f^*$. Suppose the orbit through $u^{(0)}:=\concat({\alpha^{(0)}}, {\beta^{(0)}}) \in \rr^4$ converges to $\hat{u} := \concat(\hat{\alpha}, \hat{\beta}) \in \rr^4$ and let $ \hat R := [\hat{\alpha}_1, \hat{\beta}_1] \times [\hat{\alpha}_2, \hat{\beta}_2] \subset X$ denote the equilibrium bounds guaranteed by Theorem \ref{thm:bootstrap_eqbounds}. Then exactly one of the following is true. 
	\begin{enumerate}
		\item $\hat{R}$ is a degenerate rectangle (i.e.~for $i = 1,2$, $\hat{\alpha}_i = \hat{\beta}_i$) and $f^*$ has a unique equilibrium $\hat{x} = (\hat{\alpha}_1, \hat{\alpha}_2)$ which is stable. 
		\item $\hat {R}$ is non-degenerate and $f^*$ has at least two stable equilibria. Specifically, the corners of $\hat{R}$ with coordinates
		\[
		\hat{x}_1 = (\hat{\alpha}_1, \hat{\beta}_2), \qquad \hat{x}_2 = (\hat{\beta}_1, \hat{\alpha}_2)
		\]
	 are stable equilibria of $f^*$. 
	\end{enumerate}
\end{theorem}

\begin{proof}
Define $\hat{x}_1 := (\hat{\alpha}_1, \hat{\beta}_2),\, \hat{x}_2 := (\hat{\beta}_1, \hat{\alpha}_2)$, and observe that since $(\hat{\alpha}, \hat{\beta})$ is a fixed point of $\Phi$ we have by direct computation
\begin{align*}
    H_{1,2}(\hat{\beta}_2) & = \hat{\alpha}_1 \\
    H_{2,1}(\hat{\beta}_1) & = \gamma_2^* \hat{\alpha}_2 \\
    H_{1,2}(\hat{\alpha}_2) & = \hat{\beta}_1 \\
    H_{2,1}(\hat{\alpha}_1) & = \gamma_2^* \hat{\beta}_2.
\end{align*}
It follows that 
\[
f^*(\hat{x}_1) = f^*(\hat{\alpha_1}, \hat{\beta}_2) = 
\begin{pmatrix}
\hat{\alpha}_1 - H_{1,2}(\hat{\beta}_2) \\
\gamma_2^* \hat{\beta}_2 - H_{2,1}(\hat{\alpha}_1)
\end{pmatrix}
= 
\begin{pmatrix}
0 \\
0
\end{pmatrix}
\]
\[
f^*(\hat{x}_2) = f^*(\hat{\beta}_1, \hat{\alpha}_2) = 
\begin{pmatrix}
\hat{\beta}_1 - H_{1,2}(\hat{\alpha}_2) \\
\gamma_2^* \hat{\alpha}_2 - H_{2,1}(\hat{\beta}_1)
\end{pmatrix}
= 
\begin{pmatrix}
0 \\
0
\end{pmatrix}
\]
so that $\hat{x}_1, \hat{x}_2$ are equilibria for $f^*$. Evidently, if $\hat{R}$ is degenerate then $\hat{x}_1 = \hat{x}_2$ and by Theorem \ref{thm:bootstrap_eqbounds} it follows that this is the unique equilibrium for $f^*$. On the other hand if $\hat{R}$ is nondegenerate then $\hat{x}_1$ and $\hat{x}_2$ are distinct equilibria of $f^*$. 

Next, we prove that $\hat{x}_1$ is stable as the argument for stability of $\hat{x}_2$ is similar. Observe that the derivative of $f^*$ at $\hat{x}_1$ is given by the formula
\[
Df^*(\hat{x}_1) = 
\begin{pmatrix}
-1 & H'_{1,2}(\hat{\beta}_2) \\
H'_{2,1}(\hat{\alpha}_1) & -\gamma^*_2
\end{pmatrix}
\]
which has eigenvalues satisfying 
\begin{equation}
\label{eq:corner_eigenvalues}
    z^2 + (1 + \gamma^*_2)z + \gamma^*_2 - H'_{1,2}(\hat{\beta}_2) H'_{2,1}(\hat{\alpha}_1) = 0.
\end{equation}
Since $H_{1,2}$ and $H_{2,1}$ are monotonically decreasing, the discriminant of the polynomial in Equation \eqref{eq:corner_eigenvalues} is 
\[
(1 + \gamma_2^*)^2 - 4 \paren*{\gamma_2^* - H_{1,2}'(\hat{\beta_2})H_{2,1}'(\hat{\alpha_1})} = (1 - \gamma_2^*)^2 + H_{1,2}'(\hat{\beta_2})H_{2,1}'(\hat{\alpha_1}) > 0.
\]
Hence, we can deduce that the eigenvalues are real and distinct and thus $\hat{x}_1$ is hyperbolic.

Let $f_1, f_2$ denote the two components of $f^*$. We consider points $x = (x_1, x_2)$ near the equilibrium at $\hat{x}_1$ lying on the lines defined by $x_1 = \hat{\alpha_1}$ and $x_2 = \hat{\beta_2}$. Specifically, we have the four cases: 
\begin{enumerate}
    \item If $x_1 < \hat{\alpha}_1$, then $f_1(x_1, \hat{\beta}_2) = -x_1 + H_{1,2}(\hat{\beta}_2) > -\hat{\alpha_1} + H_{1,2} (\hat{\beta}_2) = 0$.
    \item If $x_1 > \hat{\alpha}_1$, then $f_1(x_1, \hat{\beta}_2) = -x_1 + H_{1,2}(\hat{\beta}_2) < -\hat{\alpha_1} + H_{1,2} (\hat{\beta}_2) = 0$.
    \item If $x_2 < \hat{\beta}_2$, then $f_2(\hat{\alpha}_1, x_2) = -\gamma_2^*x_2 + H_{2,1}(\hat{\alpha}_1) > -\gamma_2^*\hat{\beta}_2 + H_{2,1}(\hat{\alpha}_1) = 0$.
     \item If $x_2 > \hat{\beta}_2$, then $f_2(\hat{\alpha}_1, x_2) = -\gamma_2^*x_2 + H_{2,1}(\hat{\alpha}_1) < -\gamma_2^*\hat{\beta}_2 + H_{2,1}(\hat{\alpha}_1) = 0$.
\end{enumerate}

Consequently, in a small enough neighborhood of $\hat x_1$  the flow along each axis points towards the equilibrium proving that the equilibrium is stable. 

\end{proof}

\subsection{Equilibria and saddle node bifurcations in Toggle Switch} \label{sec:saddlenode_TS}

In this Section, using the boxier box algorithm and the intuition highlighted in Remark \ref{rem:saddlenode_alg}, we want to explore the appearance of saddle nodes in the different regions of the Toggle Switch. In particular, we want to focus here on vertical paths, as introduced in Section \ref{sec:equilibria}, thus only varying the Hill coefficient for an otherwise fixed parameter. 

Interpreting the results provided by DSGRN and presented in Table \ref{tab:parameter_regions}, we can partition the parameter space into two groups: parameter region 5, $R(5)$, that is bistable, and the rest of the regions around $R(5)$, that we can call \emph{the doughnut} and is monostable.

In this Section, our aim is double: on one hand, we want to show how large Hill coefficients generate the same dynamics as the combinatorial dynamics studied by DSGRN, and on the other we posit that saddle nodes are the main driver for new equilibria in the Toggle Switch.

For these two goals, we introduce the concepts of coherency and coherency rate.

\begin{definition}
    A parameter $p\in \Lambda$ is said to be \emph{coherent} for a given Hill system if the number of stable equilibria numerically found is the same as the one determined by DSGRN for $\pi(p)$.

    The \emph{coherency rate} of a set of parameters is the percentage of coherent parameters in the set for the given Hill system.
\end{definition}

Let us underline here that, for an infinitely large Hill coefficient $d \rightarrow \infty$, all parameters are theoretically coherent since they converge to the parameters DSGRN actually computes on.

We now need a probability distribution to generate our parameters. As can be seen from Table~\ref{tab:parameter_regions} each region $R(k)$ is unbounded. With this in mind we make use of a Gaussian (as opposed to uniform) distribution for the random sampling.
To generate random parameters for a statistical analysis, we construct an unbounded distribution $\mathcal{F}$ such that samples taken from $\mathcal{F}$ span the parameter regions $R(i)\subset \Xi^*, i = 1,\dots, 9$, of the Toggle Switch in such a way that no region is significantly over- or under-represented. Thus, we create large samples of parameters knowing that $R(5)$ will be roughly represented in the sample as often as the other parameter regions. A sample from such a random distribution is built by {\tt create\_dataset\_ToggleSwitch}.

\begin{figure}[t]
	\begin{center}
		\includegraphics[width = 0.4\textwidth, keepaspectratio=true, trim={0 0 0 0},clip]{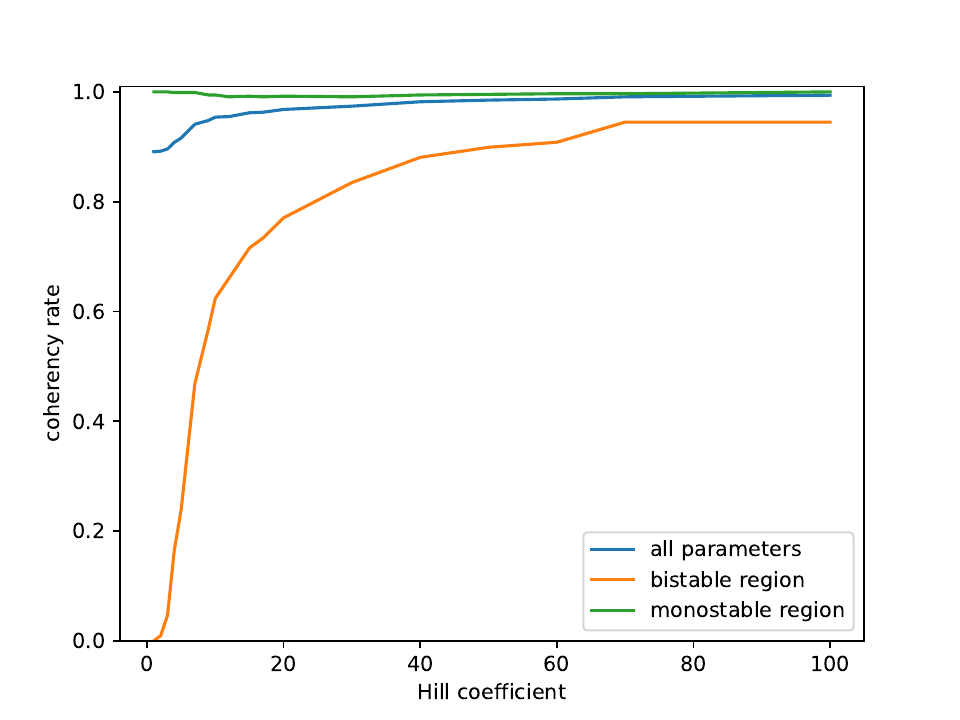}
		\caption{Coherency rate as a function of the Hill coefficient depending on the DSGRN prediction. In blue all parameters, in green only parameter predicted to be monostable, in orange parameters predicted to be bistable. As we can see, for the bistable region a choice of a large enough Hill coefficient is paramount for the detection of the bistability, while this does not hold true in the monostable region, where monostability is reliably detected independently of the Hill coefficient.
		}
		\label{fig:coherency_TS}
\end{center}
\end{figure}

In Figure \ref{fig:coherency_TS}, created by {\tt TS\_coherency}, we can grasp an understanding of the presence of equilibria depending on the region and on the Hill coefficient for the Toggle Switch. We can also compute the presence of saddle nodes. As the Hill exponent varies between 1 and 100, 93.7\% of parameters in $R(5)$ undergo a saddle node bifurcation as compared to 0.4\% of parameters outside of $R(5)$. 

\subsection{Multiple saddle-node bifurcations, hysteresis, and isolas}
\label{sec:hysteresis and isolas}

Despite the remarkable accuracy of DSGRN predictions demonstrated in Section \ref{sec:saddlenode_TS}, there exist relatively small sub-regions of parameters for which there is evidence of richer dynamics than predicted. This occurs even in simple examples such as the Toggle Switch and in this section we demonstrate our ability to find such instances and explore some of the dynamics exhibited in these instances. 

\begin{figure}[t]
	\begin{center}
		\includegraphics[width = \textwidth, keepaspectratio=true, trim={0 0 0 1.2cm},clip]{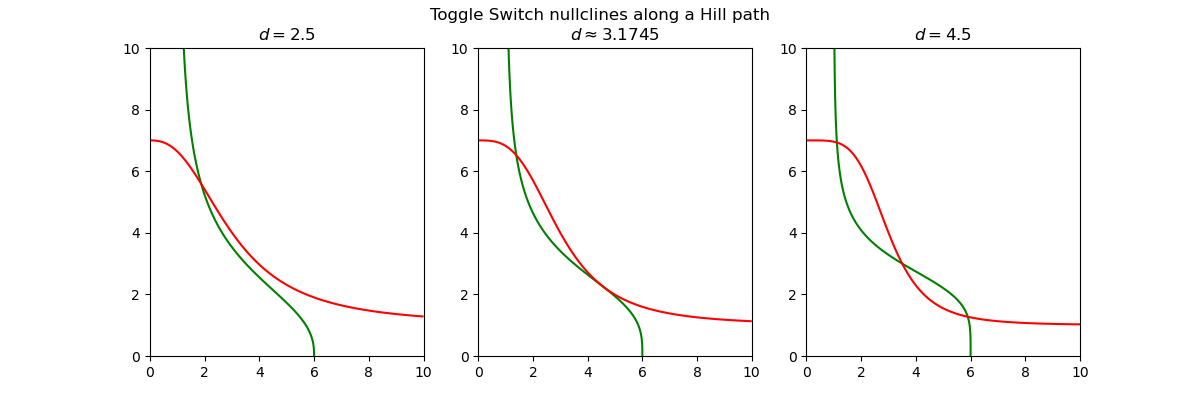}
		\caption{Toggle Switch nullclines for a parameter in $R^*(5)$ with Hill exponent varying through a saddle-node bifurcation. As both Hill functions become steeper, the nullclines tend to be closer to one another until the saddle-node bifurcation, at which point the nullclines intersect tangentially. Further increasing $d$ increases the area bounded between the nullclines and the three equilibria appear to move further apart from each other. 
		}
		\label{fig:nullclines}
\end{center}
\end{figure}

In Section \ref{sec:saddlenode_TS} our numerical analysis hints to the fact that there is a strong correlation between parameter regions and saddle node behaviour, in particular one of these two possibilities is likely to hold:
\begin{itemize}
    \item $\pi_{\Xi^*}(\xi^*) \in R^*(5)$, then along the path \textbf{r} there exists a saddle-node bifurcation for some $d\in (1,\infty)$.
    \item $\pi_{\Xi^*}(\xi^*) \notin R^*(5)$, then for all $d\in[1,\infty)$ along \textbf{r} there exists a unique equilibrium.
\end{itemize}
This second element is also supported by the following intuition. For $d > 2$ both Hill functions are sigmoidal with $d$ controlling the incline of the steepest part of the sigmoid. Therefore, for the Toggle Switch with identified Hill coefficients, it seems reasonable to guess that the nullclines of $f^*$ move ``toward each other'' monotonically as $d \to \infty$. This intuition is shown in Figure \ref{fig:nullclines} and suggests that $f^*$ may undergo at most one saddle-node bifurcation with respect to $d$.

\begin{figure}[th]
	\begin{center}
		\includegraphics[width = 0.45\textwidth]{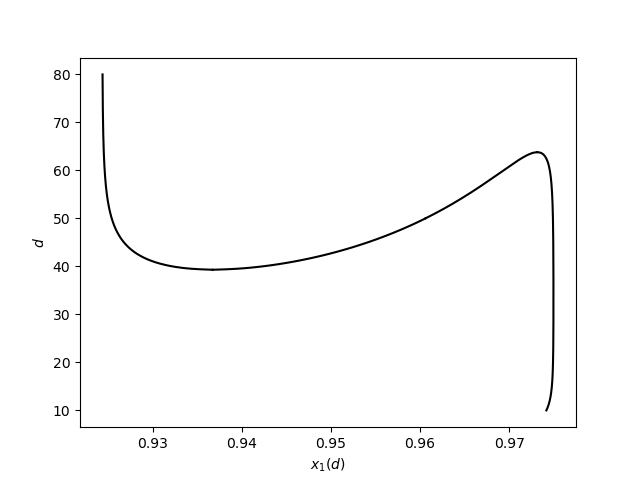}
		\includegraphics[width = 0.45\textwidth]{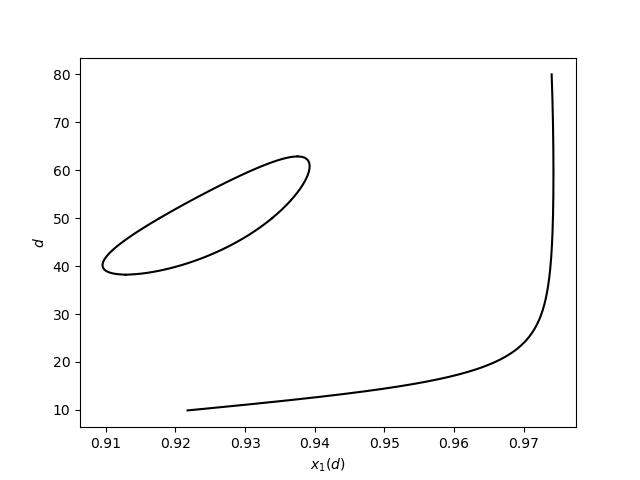}
		\caption{
			Branches of equilibria computed for the reduced Toggle Switch for varying $d$ via pseudo-arc length continuation depicting the two mechanisms which cause multiple saddle-node bifurcations. (Left) Hysteretic behaviour for the Toggle Switch found at parameter value $\lambda^*=[0.9243, 0.0506, 0.8125, 0.0779, 0.8161] $. (Right) An isola is found at $\lambda^* = [0.6470, 0.3279, 0.9445, 0.5301, 0.3908]$ indicating that multiple saddle-node bifurcations are possible without hysteretic switching.}
		\label{fig:numerical_isolas}
	\end{center}
\end{figure}

However, the techniques described in this paper provide the means to efficiently test this conjecture. Despite the fact that the parameter space is ``high'' dimensional, the combinatorial analysis suggests that if a counterexample exists, one might look for it near the boundary of $R^*(5)$. Then, using the efficient numerical methods described in this work, we easily find parameters which numerically demonstrate multiple bifurcations and surprisingly this conjecture appears to be false. 
Indeed, the numerical investigation suggests that parameters outside of $R^*(5)$ may undergo (at least) two saddle-node bifurcations along a Hill path at some Hill coefficients $1 < d_1 < d_2 < \infty$. For Hill coefficients in $[1, d_1)$ and $(d_2, \infty)$, the system is monostable and it is bistable in $(d_1 , d_2)$. This behavior is shown for two parameters in Figure \ref{fig:numerical_isolas}. 
Moreover, the transition between monostability and bistability need not be hysteretic with respect to the Hill coefficient.
Indeed, pairs of saddle-node bifurcations occurring along a single Hill path may be hysteretic (i.e.~connected by a single equilibrium branch) or form isolas which are disconnected from a ``main'' branch associated to stable equilibrium which persists for all Hill coefficients. The combination of these dynamic behaviors in this model suggests one possible cause for failures in the construction of synthetic biological switches. 
\if false
This means that, for any parameter outside $R(5)$ for which we were able to find one saddle node, we should be able to find a second one, where the two additional branches of equilibria disappear.
Considering the equilibrium branches involved in these two saddle node bifurcations, the two possible behaviour of the equilibria are either to be hysteretic, or to define an isola. The graphical difference between the two behaviours can be seen in the Figure \ref{fig:isolas}.

\begin{figure}[h]
	\begin{center}
		\includegraphics[width = 0.7\textwidth, trim= 1.5cm 14cm 4.25cm 1.7cm, clip]{sketch_isolas}
		\caption{On the left, expected behavior of a saddle-node in the center region $\mathcal{R}$. On the right, an isola,  expected behavior of equilibria when a saddle-node is detected outside the center region. Vertical axes in this case is $\hill$, while the saddle-nodes are indicated by a red diamond. Unstable equilibria are drawn as dotted lines.}
		\label{fig:isolas}
	\end{center}
\end{figure}
To the best of our knowledge, only hysteretic behaviour is expected to be found in the toggle switch \cite{kuznetsov2004synchrony}, but we were able to identify parameters showcasing both behaviours.
\fi 
In Figure \ref{fig:numerical_isolas}, we present a plot of the Hill coordinate with respect to the $x_1$ coordinate of numerically found equilibria for two  different parameters outside of the center region. We can see how they differ in the fact that, in the first parameter, a continuous function can be plotted, while in the second one, the Hill coefficient is not a function of $x_1$ and creates an isola. In both cases, the saddles have been first found numerically, using the method presented in Section \ref{sec:equilibria}, then the equilibria have been numerically continued to produce the figures.

\subsection{Degenerate saddle-node bifurcations}
\label{sec:degenerate saddle-node bifurcations}
In this section we numerically investigate and compare the saddle-node bifurcations in the Toggle Switch occurring near the ``corners'' of $R(5)$. We note how the description of the attractors given by DSGRN includes not only a count of the stable equilibria but also their relative position in a simplified phase space. This richer information can be used to further investigate the dynamical behaviour supported by a system by applying our combined combinatorial-numerical approach. 

Consider the simplest situation that all of the attractors in the Morse decomposition identified for each region of the combinatorial parameter space in Table \ref{tab:parameter_regions} is associated with a single stable equilibrium for some smooth vector field $f : (0, \infty)^2 \to \rr^2$ located in the indicated quadrant of the state space. We consider first a path through the continuous parameter space which visits the parameter regions in the order $R(5) \mapsto R(8) \mapsto R(7) \mapsto R(4) \mapsto R(5)$. In the combinatorial system, a continuous branch of equilibria can not cross the hyperplanes defined by the equations $x_1 = \gamma_1 \theta_{2,1}$ and $x_2 = \gamma_2 \theta_{1,2}$. Therefore, if $f$ is a faithful continuous representation of the combinatorial dynamics then the most parsimonious explanation for the sequence of Morse decompositions along this path is as follows. 

For parameters in $R(5)$, $f$ has (at least) one unstable equilibrium. When passing between regions $R(5)$ and $R(8)$, this unstable equilibrium and the stable equilibrium in $(1,0)$ disappear in a saddle-node bifurcation. Similarly, when passing between regions $R(5)$ and $R(4)$ this unstable equilibrium and the stable equilibrium in $(0,1)$ disappear in a saddle-node bifurcation. Therefore, along the boundary separating $R(4)$ and $R(7)$ these two saddle-node bifurcations collide in a codimension-2 cusp bifurcation and similarly on the boundary between $R(8)$ and $R(7)$. An analogous analysis of a path visiting the regions (in order) $R(5) \mapsto R(2) \mapsto R(3) \mapsto R(6 )\mapsto R(5)$ suggests another cusp bifurcation at the ``corner'' shared by $R(5)$ and $R(3)$. 

To contrast this behavior, consider a path visiting the parameter regions $R(5) \mapsto R(8) \mapsto R(9) \mapsto R(6) \mapsto R(5)$. For all regions along this path, the unique stable equilibrium is located in quadrant $(0,1)$ which is most likely modelled by the existence of a common smooth branch of equilibria which remains in this quadrant throughout the path. On the boundary between $R(5)$ and $R(6)$, the unstable equilibrium and the stable equilibrium in $(1,0)$ disappear in a saddle-node bifurcation much like on the boundary between $R(5)$ and $R(8)$. 
A similar situation occurs for paths visiting regions $R(5) \mapsto R(2) \mapsto R(1) \mapsto R(4) \mapsto R(5)$. This behaviour is not likely compatible with a cusp bifurcation at the "corner" between $R(5)$ and $R(1)$ or $R(9)$. Thus, our conjecture is that bifurcations for parameters at the ``corners'' shared by region $R(5)$ and regions $R(1), R(9)$ and those at the corners shared by $R(5)$ and regions $R(3), R(7)$ are caused by very different mechanisms.

To investigate our conjecture, suppose once again that $f^*$ is the nondimensionalized Toggle Switch Hill model defined in section \ref{sec:reducing the number} with $6$ dimensionless parameters. For visualisation purposes, we  define a suitable dimension reduction map from from $\Lambda^*$ onto the rectangle $[0, 3]^2 \subset \rr^2$ which preserves these $9$ DSGRN parameter regions and the relative distances between parameters and the region boundaries. 

To start, suppose $\xi^* \in \Xi^*$ and define the related parameters $\setof*{a_1, b_1, a_2, b_2}$ by the formulas
\[
a_1 := \ell^*_{1,2} \qquad b_1 := \ell_{1,2}^* + \delta_{1,2}^* \qquad a_2 := \frac{\ell_{2,1}^*}{\gamma_2^*} \qquad b_2 := \frac{\ell_{2,1}^* + \delta_{2,1}^*}{\gamma_2^*}.
\]
Observe that if the reduced polynomial inequalities in Table \ref{tab:parameter_regions} are expressed in terms of these new parameters, then the $9$ DSGRN parameter regions as well as their boundaries are defined by linear manifolds. For example, region $R(5)$ is given by
\[
R(5) = \setof*{(a, b) \in \rr^2 \times \rr^2 : 0 < a_1 < 1 < b_1 \ \text{and} \ 0 < a_2 < 1 < b_2 }
\]
and the boundary separating $R(5)$ and $R(6)$ is given by 
\[
\partial R(5) \cap \partial R(6) = \setof*{(a, b) \in \rr^2 \times \rr^2 : 0 < a_1 < 1 < b_1, \ a_2 = 1, \ \text{and} \ 1 < b_2 }.
\]
Implicitly this defines a nonlinear transformation $\psi : \Xi^* \to \rr^2 \times \rr^2$ given by the formula $\psi(\xi^*) = (a, b)$ which maps each DSGRN region to a unique hyperplane in $\rr^2 \times \rr^2$.  

To complete the construction, fix positive constants $\overbar{a}_1, \overbar{a}_2>1$ and define another map $g : \rr^2 \times \rr^2$ by the formulas
\begin{equation}\label{eq:2dmap}
    g_1(a,b) = 
\begin{cases}
	b_2 & \text{if} \ b_2 \leq 1 \\
	1+\frac{ 1 - a_2}{b_2 - a_2} & \text{if} \  a_2 < 1 < b_2 \\
	2 + \frac{a_2-1}{\overbar{a}_2 - 1} & \text{if} \ 1 \leq a_2 
\end{cases}
\qquad 
g_2(a,b) = 
\begin{cases}
	b_1 & \text{if} \ b_1 \leq 1 \\
	1+\frac{ 1 - a_1}{b_1 - a_1} & \text{if} \ a_1 < 1 < b_1 \\
	2 + \frac{a_1-1}{\overbar{a}_1 - 1} & \text{if} \ 1 \leq a_1
\end{cases}
\end{equation}
Observe that the bounded subset $K_{\overbar{a}_1, \overbar{a}_2} \subset \image \psi$ defined by the formula
\[
K_{\overbar{a}_1, \overbar{a}_2} := \setof*{(a, b) \in \image \psi : \abs{a_i} \leq \overbar{a}_i, 1< \overbar{a}_i,  i = 1,2}.
\]
satisfies $g\left(K_{\overbar{a}_1, \overbar{a}_2}\right) \subset [0,3]^2$.

We use the previous constructions to visualize the relative position of parameters as follows. Given a fixed collection of parameters $\setof*{\xi^*_1, \dotsc, \xi^*_M} \subset \Xi^*$ we choose $\overbar{a}_1, \overbar{a}_2$ sufficiently large so that $\psi(\xi^*_j) \in K_{\overbar{a}_1, \overbar{a}_2}$ for $1 \leq j \leq M$. Therefore, the mapping $g \circ \psi : \Xi^* \to [0,3]^2$ satisfies the following properties. 
\begin{enumerate}
\item Each of the $9$ parameter regions is mapped to a distinct unit square in $[0,3]^2$ and their relative positions are preserved. In fact, $g \circ \psi$ has been constructed so that the images of these regions are simply obtained by superimposing the graph in Figure \ref{fig:TS}(b) onto $[0, 3]^2$. We let $S(i) := g \circ \psi (R^*(i))$ denote these unit squares for $1 \leq i \leq 9$. 
\item Each boundary separating a pair of parameter regions is mapped into a line of the form $g_i = j$ with  $i \in \setof*{1,2}$ and $j \in \setof*{1,2,3}$. In other words, boundaries of parameter regions are mapped onto boundaries of the corresponding unit square and relative positions of the boundaries are also preserved. 
\item Relative proximity to boundaries is preserved. Specifically, suppose $\xi_1, \xi_2 \in R^*(j)$ and $\operatorname{dist}(\xi_1, \partial R(j)) < \operatorname{dist}(\xi_2, \partial R(j))$, then $\operatorname{dist}(g \circ \psi ( \xi_1), \partial S(j))  < \operatorname{dist}(g \circ \psi ( \xi_2), \partial S(j))$.
\end{enumerate}
Given this projection of the parameter space into a two-dimensional Euclidean space we now visually investigate the behavior of saddle-node bifurcations along the boundary of $R^*(5)$. 

For each parameter sample $\lambda \in \Xi^*$, we are interested in finding any saddle node bifurcation that happens along the path $\curve(s)$ as presented in Section \ref{sec:equilibria}.

\begin{figure}
\begin{center}
\includegraphics[width = 0.475\textwidth%, trim = 1cm 0 0 4.9cm, clip
]{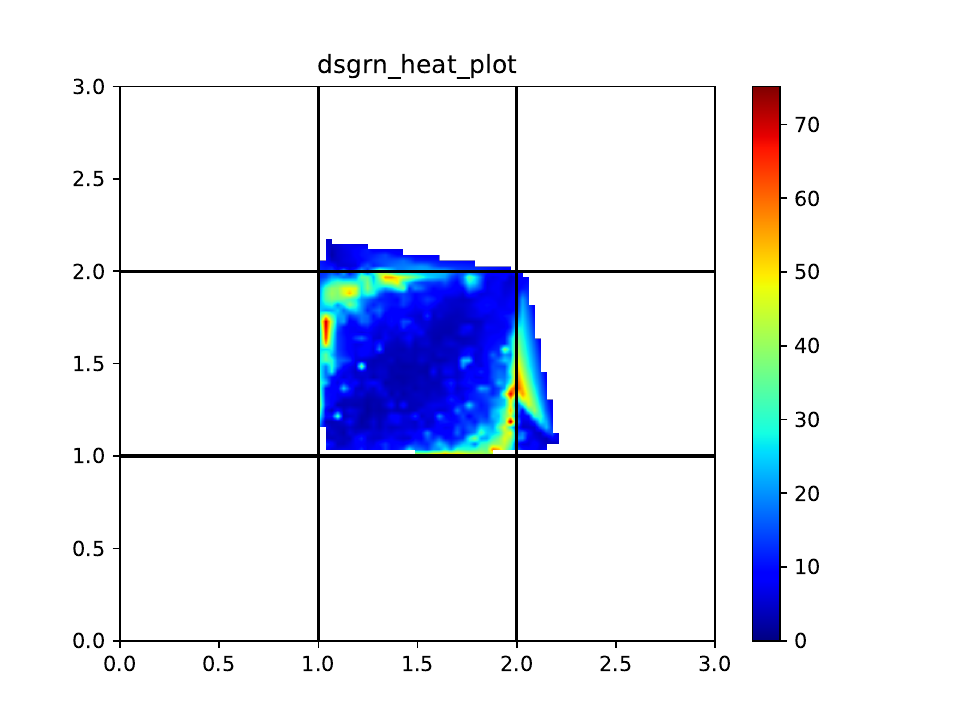}
\includegraphics[width = 0.475\textwidth%, trim = 1cm 0 0 4.9cm, clip
]{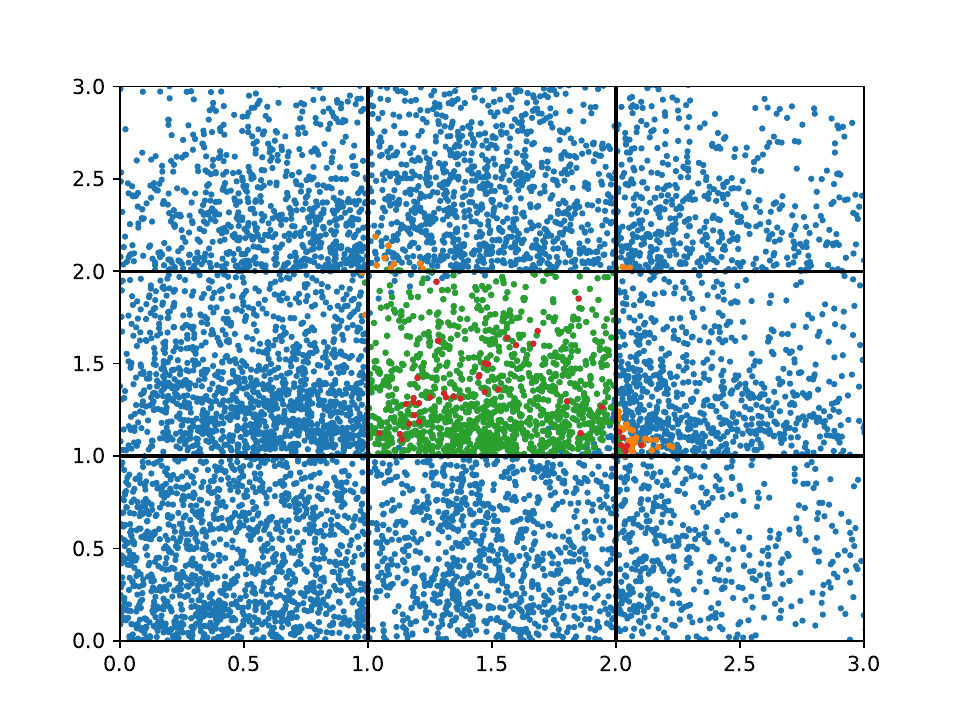}
\caption{Left: Using the projection presented in Equation \ref{eq:2dmap}, a heat map is plotted, indicating the smallest Hill coefficient undergoing a saddle node. Right: Using the same projection, parameters are plotted in blue if they don't undergo any saddle node, in green if they undergo a single saddle node, in orange if they undergo multiple saddle nodes and in red if the bisection algorithm found a saddle node that was not numerically confirmed with Equation \eqref{eq:num_saddle_node}}
\label{fig:heat_map}
\end{center}
\end{figure}

The first visual result we present is an overview of parameters that undergo saddle nodes at different Hill coefficients. For this, we find saddle nodes with respect to the Hill coefficient, then project the parameter over the $[0,3]^2$ square and represent the Hill coefficient as a heat map.
In Figure \ref{fig:heat_map}, generated with {\texttt{main\_heat\_map}}, on the left the projection \eqref{eq:2dmap} is used, while the color indicates the lowest Hill coefficient $d$ for which we could find a saddle node at the given parameter value. On the right Figure, there is a scatter plot of all parameters considered. 
These includes information on parameters that never undergo a saddle node, but also on parameters that undergo more than one saddle node, such as discussed in Section \ref{sec:hysteresis and isolas}. 

At first glance, Figure \ref{fig:heat_map} further supports the idea that choosing parameters in $R^*(5)$ has the highest likelihood of giving us a saddle node for relatively low Hill coefficient. 

Looking at this map, we observe that the bottom left of the center region seems to be the best location for a practical bistable switch, since most parameters in that section undergo a saddle node at low Hill coefficient, and thus support bistability for most large Hill coefficients. We notice how there are saddle nodes taking place outside the center region, but most of them undergo two saddle nodes, thus for high enough Hill coefficient they again present a unique stable fixed point, as discussed in Section \ref{sec:hysteresis and isolas}.

Crucially, this plot supports the claim that saddle-node bifurcations disappear into a cusp (or other codimension-2) bifurcation along a loop around the points $(1,1)$ or $(2,2)$. However, along a loop around $(1,2)$ or $(2,1)$, saddle-node bifurcations seem to disappear ``at infinity''. This is consistent with the observations we obtained from only the combinatorial dynamics using DSGRN, and that we presented at the beginning of this Section.

While the presented description encompasses the majority of the parameter values, there are still some outliers.
Indeed, we are able to find a handful of parameters in $R^*(3)$ that undergo two bifurcations. These parameters are close to the corner with $R^*(5)$ and can be seen in Figure \ref{fig:heat_map}.
These results underline that even for simple nonlinear systems there is sufficient complexity that any simple broadly sweeping overview will have points of failure. 
Nonetheless this paper presents methods to look at the general behaviour of systems, giving significant support to Figure \ref{fig:strategy}.
Let us remind the reader that each ``corner'' considered in Figure \ref{fig:heat_map} corresponds to a 4-dimensional surface in the full parameter space.

\section{The epithelial-mesenchymal transition network}
\label{sec:EMT_dyn}

Taking into account the activating $\to$ and repressing $\dashv$ edges in Figure~\ref{fig:EMT_model} the Hill model for the epithelial-mesenchymal transition network has the following nonlinearity
\begin{equation}\label{eq:EMT_Hill}
	\cH(x) := 
	\begin{pmatrix}
		H_{1,2}^-(x_2) H_{1,4}^-(x_4) \\
		H_{2,3}^-(x_3) H_{2,5}^-(x_5) \\
		H_{3,1}^+(x_1) H_{3,6}^-(x_6) \\
		H_{4,5}^-(x_5) \\
		H_{5,2}^-(x_2) H_{5,3}^+(x_3) H_{5,4}^-(x_4) \\
		H_{6,3}^-(x_3) H_{6,5}^-(x_5)
	\end{pmatrix}
\end{equation}
with the notation previously introduced. 
Since each Hill function contributes four positive parameters and each node contributes one parameter the associated parameter spaces are $\Lambda = (0, \infty)^{54}$ and $\Lambda_\bullet = (0,\infty)^{43}$.

The DSGRN parameter space does not include $d$ and thus is $(0,\infty)^{42}$.
The DSGRN parameter graph has exactly $10,368,000,000$ nodes.
Using DSGRN to compute the Morse graph for each parameter node is feasible, but carrying out significant in depth numerical exploration of such a large space is numerically unachievable. For this reason, in this section we explore the dynamical behaviour of equilibria in a broad variety of parameter regions.

In each of these parameter regions, DSGRN can provide us with a description of the dynamics at $d\rightarrow\infty$ in the form of a Morse graph. Thus, each parameter region is associated with an integer indicating the number of stable equilibria of the combinatorial dynamics.

Following Section \ref{sec:saddlenode_TS}, also for the EMT network we want to explore the changes in coherency for $d < \infty$.

In the first experiment we randomly select $N$ DSGRN parameter regions. For each region we query DSGRN for the number of stable equilibria expected at $d \rightarrow \infty$. Then we choose a single parameter in the region and compute its coherency for various choices Hill coefficients $1 < d < \infty$. This computation is carried out in the script \texttt{EMT\_DSGRNcoherency.py} with $N$ and $d$ as parameters which can be varied. 

We carried out this computation with $N = 400$ and $d = 10$ and found a coherency rate of $0.85$. Our interpretation of this computation is that even though DSGRN computes dynamics at $d \rightarrow\infty$, those computations remain relevant for large subsets of parameters even at biologically relevant Hill coefficients. 

\begin{figure}[t]
	\begin{center}
		\includegraphics[width = 0.6\textwidth, keepaspectratio=true, trim={0 0 0 0},clip]{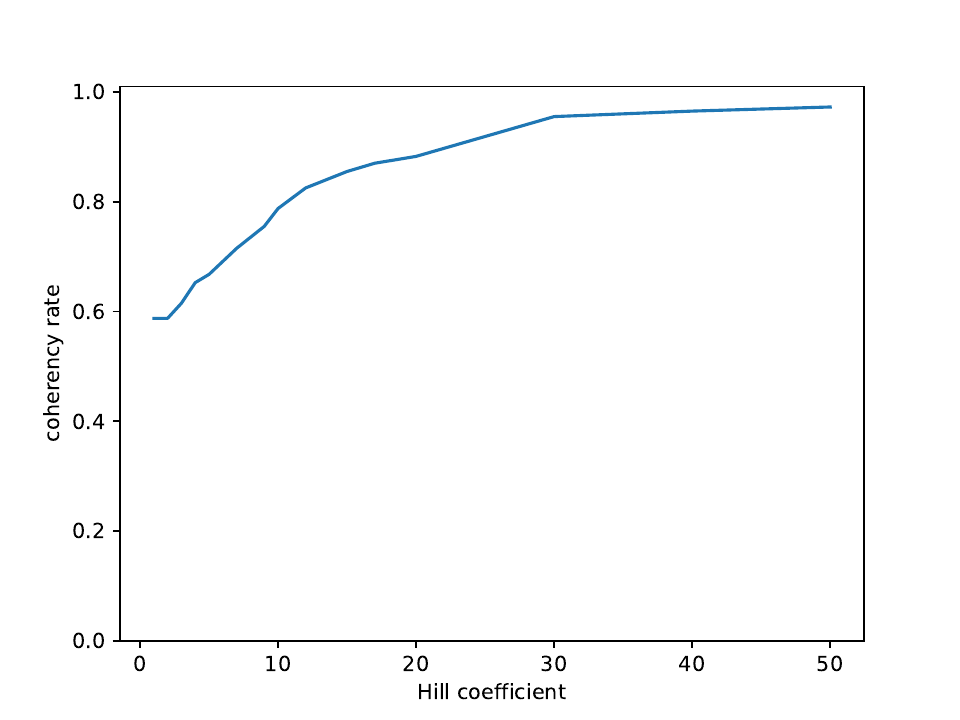}
		\caption{Coherency rate as a function of the Hill coefficient, one point per parameter region, $N = 400$.
		}
		\label{fig:coherencyVShill_singlepoint}
\end{center}
\end{figure}

Repeating this experiment for a range of Hill coefficients gives us Figure \ref{fig:coherencyVShill_singlepoint}. In this Figure, we see how, as expected, increasing the Hill coefficient increases the coherency rate. Let us remark here that for monostable regions, the coherency rate is expected to be very high starting at very low Hill coefficients, thus we now want to explore the dependency of the coherency rate on the number of stable equilibria themselves.

\begin{figure}
\centering
\begin{minipage}{.45\textwidth}
  \centering
  \includegraphics[width = 0.95\textwidth, keepaspectratio=true, trim={0 0 0 0},clip]{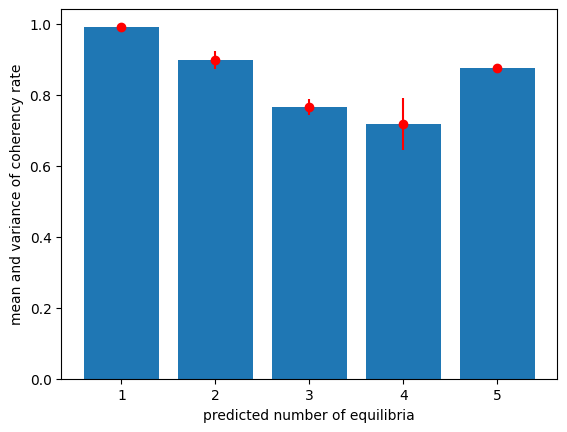}
		\caption{Mean and variance of the coherency rate depending on the number of equilibria for extremely high Hill coefficient ($d = 100$).
		}
		\label{fig:coherencyVSneqs_hill100}
\end{minipage}%
\hspace{0.05\textwidth}
\begin{minipage}{.45\textwidth}
  \centering
  \includegraphics[width = 0.95\textwidth, keepaspectratio=true, trim={0 0 0 0},clip]{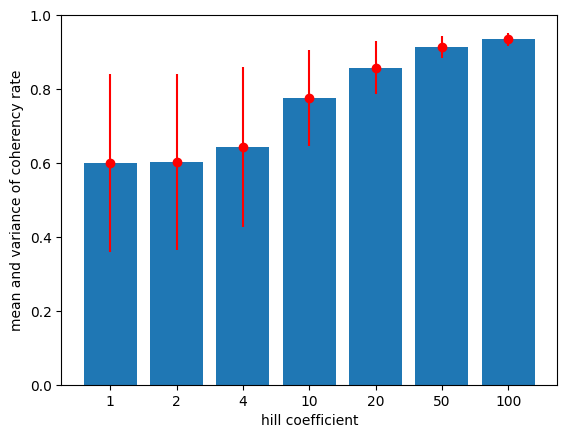}
		\caption{Mean and variance of the coherency rate depending on the Hill coefficient. It clearly shows how the dynamics for high Hill converges to the combinatorial dynamics computed by DSGRN.
		}
		\label{fig:coherencyVShill}
\end{minipage}
\end{figure}

In this second test in \texttt{EMT\_DSGRN\_nonlocal\_coherency.py}, we choose $x$ regions and for each of these we create $y$ random points in the region. Then, we compute the coherency rate of the $y$ points. Once these results computed, we can plot them in a variety of ways. First, fixing the Hill coefficient $d=100$, we check the coherency rate with respect to the number of equilibria that DSGRN detects, Figure \ref{fig:coherencyVSneqs_hill100}. Here we see that, even for an extremely high Hill coefficient, numerical limitations influence the number of found stable equilibria.

Next, we look at the influence of the Hill coefficient itself on the coherency rate. 
In Figure \ref{fig:coherencyVShill}, the coherency rate is plotted with respect to the chosen Hill coefficient. 
Here we can appreciate how the higher the Hill coefficient the higher the coherency for all parameters. 

\begin{figure}[t]
	\begin{center}
		\includegraphics[width = 0.6\textwidth, keepaspectratio=true, trim={0 0 0 0},clip]{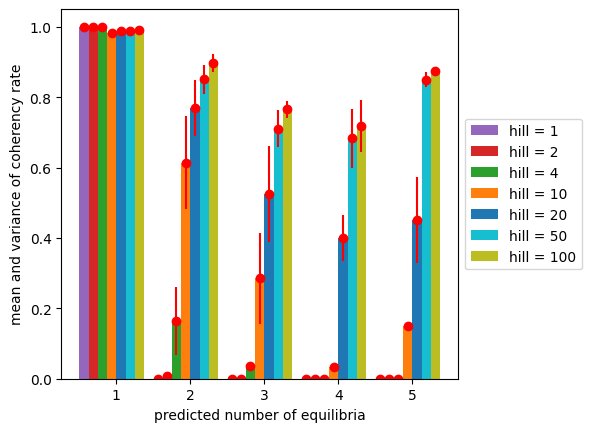}
		\caption{Mean and variance of the coherency rate depending on the Hill coefficient and the number of equilibria of the combinatorial dynamics.
		}
		\label{fig:coherencyVSneqsANDhill}
\end{center}
\end{figure}

Finally, we plot the combined results, by showing the behaviour of the coherency rate depending on the Hill coefficient and the combinatorial behaviour of the selected parameter. 
In Figure \ref{fig:coherencyVSneqsANDhill} such a histogram is presented. 
Here we see that for trivial parameters with a unique fixed point the Hill coefficient has very limited importance, while in all other cases to identify all the stable equilibria predicted by the combinatorial approach requires the use of higher Hill coefficients.

For a final in depth dynamical search, we now concentrate ourselves on only two regions, in \texttt{EMT\_pairwise\_comparison.py}. These regions are chosen such that one is defined by DSGRN as a monostable region $A$ (unique stable fixed point) and a bistable region $B$ (2 stable fixed points). Furthermore, these two regions are chosen in such a way that the monostable region has many neighboring monostable regions and the bistable has many neighboring bistable regions. This is to ensure that we are looking at a ``frontier'' between two dynamical behaviours. The chosen parameter regions for our study are $A = 14643097617$ and $B =  14643061617$. 
We generate a cloud of $200$ parameters in each of them. Their coherency rate for Hill coefficient $d = 20$ are respectively $0.975$ and $0.945$. As expected, the coherency rate for region $A$ is higher thanks to its monostability but having chosen such a high Hill coefficient the difference is small.

For each of these parameters we search for saddle node bifurcations along vertical paths. We then find that region $A$ has a $2.0\%$ prevalence of saddles, while region $B$ has $70.5\%$ prevalence of saddles. The first piece of information tells us that some saddle nodes, and thus temporary bistability, happen in monsotable regions, likely close to the border towards bistability. The second deduction is of interest if this data compared to the coherency rate. Indeed a higher percentage of parameters develops bistability for high Hill coefficient than the ones found to undergo a saddle node bifurcation. This points to a combination of numerical instability (saddle node bifurcations happen but are not detected) but also the to fact that saddle node bifurcations  are not the only source of bistability, so other more complicated bifurcations could be the generators of additional fixed points.

As a final test on our two regions $A$ and $B$, we search for horizontal saddle nodes. Since $A$ and $B$ are neighboring regions, in the DSGRN setting a unique inequality has been changed between the two of them. For our specific case, in region $A$ it holds 
$$
\theta_{4,5} < \ell_{5,2}\ell_{5,3}\ell_{5,4} 
$$
while in region $B$ we have
$$
\ell_{5,2}\ell_{5,3}\ell_{5,4}  < \theta_{4,5} \text{ and } 
\theta_{4,5} < (\ell_{5,2}+ \delta_{5,2})\ell_{5,3}\ell_{5,4}
$$
we can thus choose a horizontal path by increasing $\theta_{4,5}$ thereby crossing from region $A$ to region $B$. 

In doing so, we can look for saddle node bifurcations happening with respect to $\theta_{4,5}$, that we call horizontal saddles. In this situation, with regions $A$ and $B$ randomly chosen, we can find 200 out of 200 parameters tested for Hill coefficient 20. This indicated how horizontal saddles are the driving mechanism for the appearance of bistability.

\section{Conclusion}
\label{sec:conclusion}

The results of this paper provide additional insights into the dynamics of the Toggle Switch and EMT.
However, our goals discussed below, which we believe we have achieved, are much broader in that we have introduced novel approaches to the analysis of systems of ODEs with high dimensional parameter spaces. 
\begin{description}
\item[Integrating combinatorial and numerical approaches:] In the dynamical system community, many numerical methods have been developed for the study of dynamical systems with few parameters, but these methods are hard to apply when confronted with broader parameter spaces. For these situations, an initial overview can be gathered from the results presented by combinatorial dynamics. In this paper, we mainly consider the output of the DSGRN software, that provides us with plentiful information for our needs - mostly the number of equilibrium points depending on the parameter region, but future work could involve periodic orbits and transient dynamics. The combinatorial approach also further informs the choice of interesting \emph{paths}. In this paper we concentrated ourselves on vertical and horizontal paths, mostly driven by the simplicity of the explanation and of the code, but also because those sort of paths already cover all the possible connecting paths between combinatorial regions. This being said, we expect interesting developments in future years on the open problem of the optimal choice of path, where applications will be a driving force in this regard.

\item[The boxy box algorithm:] this novel algorithm (Section \ref{sec:the bootstrap algorithm} and Section \ref{sec:numerical_analysisTS}) allows us to bound all the equilibria of a broad class of dynamical systems. In our use case, this algorithm is fundamental to slash the computational complexity of our calculations, allowing us to straightforwardly detect monostability and usually significantly limit the use of other  equilibrium searches algorithms, thus decreasing the impact of the \emph{curse of dimensionality} in the study of large systems.

\item[Non-rigorous radii polynomial approach:] on one hand we are confronted with the problem of finding all equilibria, while on the other we need to solve the problem of recognising two numerical approximations of the same analytical equilibrium. In this paper, we tackle this problem with a non-rigorous application of the radii polynomial approach (Section \ref{sec:isolating equilibria}), a method first applied in the validated numerics community to bound computational error. While the rigorous computation of the error bound is not in the scope of our interests, the application of this second algorithm makes our results more trustworthy, while solving the problem of the numerical doubling of equilibria.
\end{description}

A reason for the complexity of the treated material is the stiffness of Hill systems for high Hill coefficients. This difficulty has required a certain dose of ingenuity, but is also the key element for the construction of the correlation between Hill models and combinatorial models. 

The combination of these varied tools allows not only for the generation of new and exciting results in the bifurcation and dynamical analysis of Hill systems, but also for an incredibly efficient code. Specifically, most of our figures have been generated in a matter of a handful of minutes on a laptop. In particular, any single equilibrium search problem takes a fraction of a second, being basically instantaneous when monostability is detected, while a saddle node search takes around half a second. As reference, the most computationally intensive pictures, Figure \ref{fig:coherencyVSneqsANDhill} and Figure \ref{fig:heat_map} took 2 and 6 hours respectively on a MacBook Pro 2019 2.3 GHz 8-Core Intel Core i9.

\section{Acknowledgements}

KM was partially supported by NIH 5R01GM126555-01, the National Science Foundation under awards DMS-1839294 and HDR TRIPODS award CCF-1934924, DARPA contract HR0011-16-2-0033,  and Air Force Office of Scientific Research under award number FA9550-23-1-0011 and FA9550-23-1-0400. 
KM was also supported by a grant from the Simons Foundation.
EQ was partially supported by the National Science Foundation under awards NIH 5R01GM126555-01 and by the German Science Foundation (Deutsche Forschungsgemeinschaft, DFG) via the Walter-Benjamin Grant QU 579/1-1.

\bibliographystyle{unsrt}
\bibliography{hill_paper}

\end{document}